\documentclass[reqno]{amsart}
\usepackage{amssymb,amsmath,enumerate,amscd,mathrsfs}

\numberwithin{equation}{section}

\newtheorem{theorem}[equation]{Theorem}
\newtheorem{proposition}[equation]{Proposition}
\newtheorem{lemma}[equation]{Lemma}
\newtheorem{corollary}[equation]{Corollary}

\theoremstyle{definition}
\newtheorem{definition}[equation]{Definition}
\newtheorem{notation}[equation]{Notation}
\newtheorem{remark}[equation]{Remark}

\DeclareMathOperator{\sym}{ \sigma\!\!\!\sigma}

\def\A{\mathcal A}
\def\B{\mathcal B}
\def\C{\mathbb C}
\def\D{\mathscr D}
\def\F{\mathscr F}
\def\Z{\mathbb Z}
\def\K{\mathcal K}
\def\L{\mathscr L}
\def\N{\mathbb N}
\def\R{\mathbb R}
\def\S{\mathscr S}
\def\x{\mathbf x}
\def\m{\mathfrak m}

\def\cusym{\,{}^{\textup{cu}}\!\sym}
\def\cusymb{\,{}^{\textup{cu}}\!\sym_{\partial}}
\def\cuT{\,{}^{\textup{cu}}T}
\def\cupi{\,{}^{\textup{cu}}\hspace{-1.5pt}\pi}

\def\Mbar{\overline{M}}
\def\regbdry{\partial_{\textup{reg}}}
\def\singbdry{\partial_{\textup{sing}}}
\def\sing{\textup{sing}}
\def\reg{\textup{reg}}
\def\cu{\textup{cu}}
\def\Psicu{{}^{\textup{cu}}\Psi}
\def\Hcu{H_{\textup{cu}}}
\def\L2cu{L^2_{\textup{cu}}}
\def\Scu{{}^{\cu}\S}

\def\eps{\varepsilon}
\def\id{\textup{I}}
\def\open#1{{\mathring{#1}}}
\def\st{;\;}

\def\cl{\textup{cl}}

\DeclareMathOperator{\Diff}{Diff}
\DeclareMathOperator{\supp}{supp}
\DeclareMathOperator{\Hom}{Hom}
\DeclareMathOperator{\codim}{codim}

\begin{document}
\title{Elliptic boundary problems on manifolds with polycylindrical ends}
\author{Thomas Krainer}
\address{Institut f\"ur Mathematik\\ Universit\"at Potsdam\\ Postfach 60 15 53\\
D-14415 Potsdam, Germany}
\email{krainer@math.uni-potsdam.de}
\begin{abstract}
We investigate general Shapiro-Lopatinsky elliptic boundary value problems on
manifolds with polycylindrical ends. This is accomplished by compactifying such
a manifold to a manifold with corners of in general higher codimension, and we then deal
with boundary value problems for cusp differential operators. We introduce
an adapted Boutet de Monvel's calculus
of pseudodifferential boundary value problems, and construct parametrices
for elliptic cusp operators within this calculus.
Fredholm solvability and elliptic regularity up to the boundary and up to infinity
for boundary value problems on manifolds with polycylindrical ends follows.
\end{abstract}

\subjclass[2000]{Primary: 58J32; Secondary: 35J70, 35J40, 35J55}
\keywords{Boundary value problems, manifolds with cylindrical ends, manifolds with
corners, degenerate elliptic operators, cusp calculus, parametrices}

\maketitle


\section{Introduction}
\label{sec-Introduction}

\noindent
This paper is concerned with the analysis of partial differential equations on noncompact
manifolds with boundary. More precisely, we are interested in Fredholm solvability, regularity, and
asymptotics of solutions for general elliptic boundary value problems on manifolds with polycylindrical ends and
boundary.

On smooth compact manifolds with boundary, it is known from classical results that Fredholmness and regularity of solutions
are governed by a tuple of principal symbols (cf.~\cite{AgmonDouglisNirenberg}): The homogeneous principal symbol of
the operator, which relates to the interior of the manifold, and the principal boundary symbol.
The index problem on compact manifolds with boundary led Boutet de Monvel in \cite{BdM71} to introduce an algebra of
pseudodifferential boundary value problems which contains the parametrices of elliptic operators
(see \cite{GrubbBuch, SzWiley98}). In particular, Fredholm solvability and regularity of solutions
follow as immediate consequences from the mapping properties of the operators in Boutet de Monvel's algebra -- this
represents a major philosophical viewpoint of pseudodifferential operator theory in general.
While the main analytic results for elliptic partial differential equations were already obtained
in \cite{AgmonDouglisNirenberg, AV63, Seeley}, Boutet de Monvel's calculus had a substantial impact on many
investigations in index and spectral theory on compact manifolds with boundary.

A significant amount of the literature on elliptic theory on noncompact manifolds focusses on
index or spectral theory related to complete Riemannian manifolds (see, e.g., \cite{AmLauNis, LauterMoroianu,
Loya, MazzeoMelrose, RBM2, MelroseNistor, Mueller, Vaillant}). The incomplete case still constitutes a major
challenge, and even for apparently ``simple'' cases like cone or edge singularities many central questions are
still unsolved (see, e.g., \cite{GKM1, GKM2, Lesch}). Substantial progress in singular elliptic theory for the incomplete situation has been achieved
over the past 10--15 years by Schulze and collaborators, see \cite{KaSchu03, SzNorthHolland, SzWiley98}.

It is known from all these works that on noncompact manifolds additional conditions that govern
Fredholmness, regularity, and asymptotics of solutions to elliptic equations at the noncompact ends
are to be expected. These conditions
are typically given in terms of operator families, and it makes sense to regard these families as principal symbols at infinity that are
associated with the operator. The principal boundary symbol of an elliptic boundary value problem mentioned above is
a classical example of such a family when we regard the boundary as being located at infinity. Other examples
are the conormal symbol for cone operators, the principal edge symbol for edge operators (cf.~\cite{SzNorthHolland}),
and the normal operators for totally characteristic or cusp operators (\cite{MazzeoMelrose, MM, LauterMoroianu}).

The topological concept of a manifold with corners from \cite{RBM2} turned out to be particularly well suited for many
problems on noncompact manifolds as it provides a natural habitat: The original manifold is contained as a dense submanifold
in the manifold with corners, and the behavior of the operators and of the solutions to the equations at the boundary reflects
the behavior at infinity of the original problem. Typically, there is an additional operator-valued condition that governs
ellipticity for each hypersurface of the boundary. All the examples mentioned above can be formulated in these terms.

Elliptic theory related to complete Riemannian geometry is often concerned with differential operators that are generated by the
smooth functions on a manifold with corners and a certain Lie algebra of vector fields (cf.~\cite{AmLauNis, MazzeoMelrose}),
while singular elliptic problems typically have singular coefficients. Despite of structural similarities in many
cases, this leads in general to a completely different analysis: Singular problems typically
induce extra conditions along the hypersurfaces of the boundary --- in the theory of boundary value problems these
are the boundary conditions, and pseudodifferential conditions of trace and potential type can be imposed at edges
to obtain a well posed problem for edge-degenerate operators (cf.~\cite{SzNorthHolland, KaSchu03}).

Elliptic theory on manifolds with polycylindrical ends without boundary relates to totally characteristic
or cusp operators on manifolds with corners as studied in \cite{LauterMoroianu, MazzeoMelrose, MelroseNistor}.
The situation we are interested in corresponds to operators that are totally characteristic
or cusp-degenerate only at \emph{some} boundary hypersurfaces, while they are regular at other hypersurfaces.
We will focus here on cusp operators rather than totally characteristic ones as they provide a larger reservoir of
admissible problems on manifolds with polycylindrical ends. At the hypersurfaces on which the operators are regular
elliptic we are imposing boundary conditions. Motivated by questions arising in cutting and pasting in $b$-geometry,
Loya and Park investigated in \cite{LoyaPark04} the Cauchy data space and the Calder{\'o}n projector at regular
hypersurfaces for totally characteristic operators of Dirac-type on manifolds with corners. As is the case also in
the classical setting of smooth compact manifolds, their work is not covered by our framework unless the boundary condition
for the Dirac operator under consideration is Shapiro-Lopatinsky elliptic (e.g., if it is local).

It is worthwhile emphasizing that boundary value problems on manifolds with polycylindrical ends are substantially different
from the corresponding problems on manifolds with singular geometric corners (i.e. manifolds
with corners endowed with a smooth nondegenerate metric).
This is reflected in our situation by the condition that the regular boundary hypersurfaces of the manifold with corners are not allowed
to have nontrivial intersections. In the presence of singular geometric corners, the natural
topological habitat is a more singular manifold with corners obtained by introducing polar
coordinates. In this case, some hypersurfaces are fibred, and --- the most significant drawback ---
the operators have singular coefficients.
To give a sufficiently general understanding to this situation remains a challenging open
problem.

This paper is organized as follows: In Section \ref{sec-ManifoldswithCorners} we briefly recall the definition of
manifolds with corners and describe the topological requirements for the boundary hypersurfaces. Moreover, we introduce
the class of cusp operators that are cusp-degenerate at the \emph{singular part}, and regular at the \emph{regular part}
of the boundary of a manifold with corners, and formulate the general boundary value problem for such operators.

From the point of view of the analysis of partial differential equations, the main result of this paper are the
characterization of ellipticity for boundary value problems given by the Definitions \ref{intcuspellipt} and
\ref{ellipticitybvp}, and Theorem \ref{diffbvpFredholm} which asserts that elliptic boundary value problems
for cusp operators are Fredholm in the natural scale of weighted cusp Sobolev spaces. Our theorem gives a precise statement as regards elliptic regularity and asymptotics of solutions up to the boundary and up to infinity, i.e. up to the regular and singular parts of the boundary.

We find that ellipticity is governed by the principal symbol (extended in a
suitable way up to the boundary), and the following additional operator families associated with each boundary hypersurface $H$:
\begin{itemize}
\item If $H$ is a regular hypersurface, ellipticity is governed by the (cusp-)principal boundary symbol.
\item If $H$ is singular, but intersects nontrivially with some regular hypersurface, then ellipticity is governed
by a family of boundary value problems for cusp operators on $H$, the conormal symbol associated with $H$.
\item If $H$ is singular and has empty intersection with the regular part of the boundary, then ellipticity
is governed by the conormal symbol, which in this case is a family of cusp differential
operators on $H$.
\end{itemize}
It should be noted that we recover the classical Boutet de Monvel algebra as corresponding to the case that the
singular part of the manifold with corners is empty, while Kondratyev's theory of elliptic boundary value problems
on conic manifolds or manifolds with mere cylindrical ends (cf. \cite{KaSchu03, K1, KraiParabBVP, KraConicResolvents,
MazNazPlam, MiNistor, SchroSchu}) corresponds in our framework to a particular case of manifolds with corners of codimension two.

The results mentioned are obtained by embedding the problems into a suitable calculus of pseudodifferential boundary
value problems for cusp operators, and constructing a parametrix for elliptic operators inside this calculus.
We set up the calculus in Section \ref{sec-Pseudocalculus}, which is the longest section of this paper.
As pointed out earlier, the pseudodifferential calculus is of independent interest in itself as it makes possible further
investigation in the direction of index and spectral theory for elliptic boundary value problems on manifolds with
polycylindrical ends.


\section{Manifolds with corners and elliptic boundary value problems}
\label{sec-ManifoldswithCorners}

\subsection{Manifolds with corners}\label{subsec-Manifoldswithcorners}

We briefly review the definition of manifolds with corners from \cite{RBM2}:

\begin{definition}\label{defmfdcorners}
An $n$-dimensional manifold with corners is a compact topological manifold $\Mbar$ with boundary
such that there exists a smooth $n$-dimensional manifold $\tilde{M}$ without boundary that
contains $\Mbar$, and smooth functions $\x_j : \tilde{M} \to {\mathbb R}$, $j = 1,\ldots,N$,
with $\Mbar = \{x \in \tilde{M};\; \x_j \geq 0,\; j=1,\ldots,N\}$, and on
$H_{i_1}\cap\ldots\cap H_{i_k}$, where $H_j = \{x \in \tilde{M};\; \x_j = 0\}$,
the differentials $d\x_{i_1}\wedge\ldots\wedge d\x_{i_k}$ are nonzero for all collections
 $i_1 < \ldots < i_k$. Without loss of generality we assume that $N$ is minimal, and so
$H_i \cap \Mbar \neq \emptyset$ for all $i = 1,\ldots,N$. In the sequel we will in general
use the notation $H_i$ when we just mean the hypersurface $H_i \cap \Mbar$ in the boundary
of $\Mbar$.

A point $p \in \Mbar$ is called a codimension $k$ point, $k \in \N$, if it lies in the
intersection of $k$ distinct hypersurfaces $H_{i_1}\cap \ldots \cap H_{i_k}$, and $k$ is
maximal with this property. The points in the
interior $M = \open{\Mbar}$ are by convention the points of codimension zero. The
codimension of a manifold with corners is defined as
$$
\codim\Mbar = \max\{k \in \N;\; \exists p \in \Mbar \textup{ of codimension $k$}\},
$$
while a manifold with corners of codimension zero is by convention a compact smooth manifold
without boundary. Note that the manifolds with corners of codimension one are just the smooth
compact manifolds with boundary.

By possibly changing the defining functions $\x_i$, we also assume that there exists $\eps > 0$
such that each boundary hypersurface $H_i$ of $\Mbar$ has a collar neighborhood diffeomorphic to
$[0,\eps)\times H_i \cong \Mbar$, and the defining function
$\x_i$ coincides in this neighborhood with the projection to the coordinate $x_i \in [0,\eps)$.
Throughout this paper these collar neighborhoods and the defining functions $\x_i$ are
henceforth fixed. Observe, moreover, that each boundary hypersurface is itself a manifold
with corners of codimension at most $\codim\Mbar - 1$. We will later make use of this collar
neighborhood structure together with an induction on the codimension of a manifold with
corners in order to define a pseudodifferential calculus that is adapted to our problem
at hand.

Finally note that we have in a canonical way well defined notions of $C^{\infty}$-functions,
tangent and cotangent bundle, as well as general smooth vector bundles and their sections
on $\Mbar$ simply by restriction from $\tilde{M}$.
\end{definition}

Our focus in this paper is the investigation of general Shapiro-Lopatinsky elliptic
boundary value problems on \emph{manifolds with polycylindrical ends}. Loya and Park studied
in \cite{LoyaPark04} Dirac-type operators on such noncompact configurations, and they call them
``manifolds with multi-cylindrical end boundaries'' (they require an additional but essentially
unnecessary topological condition for the hypersurface where the boundary condition is imposed).

The relation between manifolds with polycylindrical ends and manifolds with corners is that
the latter are compactifications of the beforementioned ones. More precisely, any
diffeomorphism $(0,\eps) \cong (-\infty,\eps)$ which maps zero to $-\infty$ can be used to push
the boundary hypersurfaces $H_i$ to minus infinity in view of the collar neighborhood structure
$[0,\eps)\times H_i \cong \Mbar$ near the boundary.
Hence the interior $M$ of $\Mbar$ is in a natural way a manifold with polycylindrical ends, which
is compactified to the manifold with corners $\Mbar$ by attaching the boundary hypersurfaces
$H_i$ at infinity as specified by the diffeomorphism (the notion of manifolds with
polycylindrical ends is in view of Definition \ref{defmfdcorners} now self-explanatory).

The situation that we are interested in corresponds to the case where not all boundary
hypersurfaces of $\Mbar$ are pushed to infinity, respectively stem from compactified cylindrical
ends. It makes sense to view the portion of the boundary $\partial\Mbar$ which arises
from compactification of noncompact ends as the \emph{singular part $\singbdry\Mbar$} of
$\partial\Mbar$, because in the study of elliptic operators the noncompactness is reflected
by a degeneracy on that part of the boundary. On the remaining part of $\partial\Mbar$ the
operators are nondegenerate, and we are asking for an elliptic boundary condition to be fulfilled
there. Hence this part of $\partial\Mbar$ is considered the \emph{regular part $\regbdry\Mbar$}
of the boundary.

Let us be more precise about the topological requirements:

\begin{definition}
Let $\partial\Mbar = \bigl(H_{1} \cup \ldots \cup H_{\ell}\bigr) \cup H_{\ell+1} \cup \ldots \cup H_N$, and
assume that $H_{i} \cap H_{j} = \emptyset$ for $i \neq j$, $i,j = 1,\ldots,\ell$. Then
$\regbdry\Mbar = \bigcup\limits_{i=1}^{\ell}H_{i}$ is an admissible choice of
boundary hypersurfaces where we can impose boundary conditions, and
$\singbdry\Mbar = \bigcup\limits_{i=\ell+1}^{N}H_i$ is the singular part of $\partial\Mbar$.
Observe, in particular, that $\overline{N} = \Mbar \setminus \singbdry\Mbar$ is a smooth manifold with boundary
$\partial\overline{N} = \regbdry\Mbar \setminus \singbdry\Mbar$.

Let $\x_{\reg} = \x_1\cdot\ldots\cdot \x_\ell$ be the total defining function for $\regbdry\Mbar$,
and $\x_{\sing} = \x_{\ell+1}\cdot\ldots\cdot \x_N$ be the total defining function for $\singbdry\Mbar$.
As is custom, we write $\x_{\reg}^{\alpha} = \x_{1}^{\alpha_{1}}\cdot\ldots\cdot \x_\ell^{\alpha_\ell}$
for $\alpha = (\alpha_1,\ldots,\alpha_\ell) \in \R^{\ell}$, and correspondingly so
for $\x_{\sing}$.
\end{definition}

\subsection{Cusp differential operators and Sobolev spaces}\label{subsec-cuspdifferentialoperators}

Every codimension $k$ point $p \in \Mbar$ has a coordinate neighborhood of the form
$[0,\eps)^k\times\Omega$ with local coordinates $\Omega \subset \R^{n-k}$ and, after renumbering
the $\x_i$'s, $(\x_1,\ldots,\x_k) \in [0,\eps)^k$. If $p \in \regbdry\Mbar$, then there
is only one hypersurface $H \subset \partial\Mbar$ with $p \in H \subset \regbdry\Mbar$
by assumption about the regular part of the boundary, and without loss of generality let this
hypersurface be $H_1 = \{\x_1 = 0\}$.

A \emph{cusp differential operator of order $m \in \N_0$} is a differential operator
$A \in \Diff^m(\tilde{M})$ restricted to $M$, which in coordinates near each codimension
$k$ point $p \in \Mbar$ is of the form
\begin{align}
\label{cuspsingular}
A &= \sum\limits_{\substack{(\alpha,\beta) \in \N_0^n \\ |\alpha| + |\beta| \leq m}}
a_{\alpha,\beta}(x_1,\ldots,x_k,y)(x_1^2D_{x_1})^{\alpha_1}\ldots(x_k^2D_{x_k})^{\alpha_k}D_y^{\beta} \\
\intertext{if $p \notin \regbdry\Mbar$, or}
\label{cuspregular}
A &= \sum\limits_{\substack{(\alpha,\beta) \in \N_0^n \\ |\alpha| + |\beta| \leq m}}
a_{\alpha,\beta}(x_1,\ldots,x_k,y)D_{x_1}^{\alpha_1}(x_2^2D_{x_2})^{\alpha_2}\ldots(x_k^2D_{x_k})^{\alpha_k}D_y^{\beta}
\end{align}
if $p \in \regbdry\Mbar$, where $a_{\alpha,\beta} \in C^{\infty}([0,\eps)^k\times\Omega)$.
A change of variables $t_i = -1/x_i$ in \eqref{cuspsingular} gives
\begin{equation}\label{cylendrepr}
A = \sum\limits_{\substack{(\alpha,\beta) \in \N_0^n \\ |\alpha| + |\beta| \leq m}}
a_{\alpha,\beta}(-1/{t_1},\ldots,-1/{t_k},y)D_{t_1}^{\alpha_1}\ldots D_{t_k}^{\alpha_k}D_y^{\beta},
\end{equation}
where the function $a_{\alpha,\beta}(-1/{t_1},\ldots,-1/{t_k},y)$ is a classical symbol separately
in each coordinate $t_i$ as $t_i \to -\infty$ (analogously for \eqref{cuspregular}). Thus
cusp (pseudo-)differential operators on manifolds with corners are associated
with the analysis on manifolds with polycylindrical ends by means of the diffeomorphism
$t = -1/x$. Another much more common setup are totally characteristic or $b$-operators which are
associated with the transformation $t = \log x$ (\cite{RBM2, MM, SchroSchu, SzNorthHolland}).
In the $b$-setup the derivatives $x_i^2D_{x_i}$ in \eqref{cuspsingular} and \eqref{cuspregular}
have to be replaced by $x_iD_{x_i}$. The setup of totally characteristic operators is more restrictive
in the sense that it is applicable to a strictly smaller class of operators on manifolds with
polycylindrical ends, because the coefficients $a_{\alpha,\beta}$ in the $t = \log x$ coordinates
have exponential asymptotics as $t_i \to -\infty$.
By combining both diffeomorphisms, every $b$-operator can be transformed
into a cusp operator and can be treated fully satisfactory using the cusp setup.
More about standard cusp pseudodifferential operators (the case $\partial\Mbar = \singbdry\Mbar$)
can be found in \cite{LauterMoroianu, MazzeoMelrose, MelroseNistor}.

The equations \eqref{cuspsingular} and \eqref{cuspregular} show
that the cusp vector fields (homogeneous real first order cusp differential operators)
are a finitely generated projective module over $C^{\infty}(\Mbar)$, and consequently are the
space of sections of a smooth vector bundle $\cuT\Mbar \to \Mbar$, the cusp tangent bundle,
which on $\overline{N} = \Mbar\setminus \singbdry\Mbar$ is canonically isomorphic to the tangent bundle $T\overline{N}$.
Locally near a codimension $k$ point $p \notin \regbdry\Mbar$, a frame for this bundle is induced by the vector fields
$x_j^2\partial_{x_j}$, $j = 1,\ldots k$, and $\partial_{y_i}$, $i = 1,\ldots n-k$,
or, if $p \in \regbdry\Mbar$, $\partial_{x_1}$ instead of $x_1^2\partial_{x_1}$.

Let $\cuT^*\Mbar$ be the cusp cotangent bundle, i.e. the dual of $\cuT\Mbar$. The
\emph{cusp-principal symbol} $\cusym(A)$ of a cusp differential operator $A$ is well defined
and a homogeneous function on $\cuT^*\Mbar \setminus 0$.
If $A$ is represented in local coordinates near a codimension $k$ point
according to \eqref{cuspsingular} or \eqref{cuspregular}, then the cusp-principal symbol
takes the form
\begin{equation}\label{cuspprincsymbol}
\cusym(A) \equiv \sum\limits_{\substack{(\alpha,\beta) \in \N_0^n \\ |\alpha| + |\beta| = m}}
a_{\alpha,\beta}(x_1,\ldots,x_k,y)\xi^{\alpha}\eta^{\beta}
\end{equation}
for $0 \neq (\xi,\eta) \in \R^{k}\times\R^{n-k}$ and $(x_1,\ldots,x_k,y) \in [0,\eps)^k\times\Omega$.

Moreover, every cusp differential operator $A$ has a \emph{cusp-principal boundary symbol}
$\cusymb(A)$, a family of operators
\begin{equation}\label{cuspbdrysymbol}
\cusymb(A) \in C^{\infty}(\cuT^*\regbdry\Mbar\setminus 0, \Hom(\Scu,\Scu)).
\end{equation}
Recall that each boundary hypersurface of $\partial\Mbar$ is again a manifold with corners,
and therefore the cusp cotangent bundle $\cuT^*\regbdry\Mbar \to \regbdry\Mbar$ is
well defined. $\Scu \to \cuT^*\regbdry\Mbar$ is a vector bundle with fibre
$\S(\overline{\R}_+)$, and it can be regarded as the space of rapidly decreasing functions in
the fibres of the inward pointing half of the conormal bundle of $\regbdry\Mbar$ in
$\cuT^*\Mbar$. More precisely, if $p \in \regbdry\Mbar$ and \eqref{cuspprincsymbol} is
a local representation of the cusp-principal symbol near $p$, then
\begin{equation}\label{cuspbdrysymbollocal}
\cusymb(A) \equiv \sum\limits_{\substack{(\alpha,\beta) \in \N_0^n \\ |\alpha| + |\beta| = m}}
a_{\alpha,\beta}(0,x_2,\ldots,x_k,y)D_{x_1}^{\alpha_1}\xi_2^{\alpha_2}\ldots\xi_k^{\alpha_k}\eta^{\beta} :
\S(\overline{\R}_+) \to \S(\overline{\R}_+)
\end{equation}
for $0 \neq (\xi',\eta) \in \R^{k-1}\times\R^{n-k}$, where $\xi' = (\xi_2,\ldots,\xi_k)$,
and $(x_2,\ldots,x_k,y) \in [0,\eps)^{k-1}\times\Omega$. Observe that the cusp-principal
boundary symbol is \emph{twisted} or \emph{$\kappa$-homogeneous} in the sense that
$$
\cusymb(A)(x',y,\varrho\xi',\varrho\eta) = \varrho^m\kappa_{\varrho}\cusymb(A)(x',y,\xi',\eta)\kappa_{\varrho}^{-1} :
\S(\overline{\R}_+) \to \S(\overline{\R}_+)
$$
for $\varrho > 0$, where $\kappa_{\varrho} : \S(\overline{\R}_+) \to \S(\overline{\R}_+)$
is the normalized dilation group action, i.e.
\begin{equation}\label{kappagroup}
\bigl(\kappa_{\varrho}u\bigr)(x_1) = \varrho^{1/2}u(\varrho x_1).
\end{equation}
The observation that the classical principal boundary symbol is twisted homogeneous led Schulze
to systematically study pseudodifferential operators with opera\-tor-valued symbols that obey twisted symbol
estimates, and these are nowadays widely applied in singular pseudodifferential operator theory
(see, e.g., \cite{SzNorthHolland, SzWiley98}).

It makes sense to regard the cusp-principal symbol \eqref{cuspprincsymbol} and the cusp-principal boundary symbol
\eqref{cuspbdrysymbol} as extensions of the principal symbol and the principal boundary symbol of $A$ from
$T^*\overline{N} \setminus 0$ and $T^*\partial\overline{N} \setminus 0$ to the cusp cotangent bundles.

For each hypersurface $H \subset \singbdry\Mbar$ there is an associated \emph{conormal symbol} or
\emph{normal operator} $N_H(A)(\tau)$ to $A$, which is a family of cusp differential
operators on the manifold with corners $H$ depending on the parameter $\tau \in \R$. If
$H = H_k = \{\x_k = 0\}$ in the local representation \eqref{cuspsingular} or \eqref{cuspregular} of $A$
near a point $p \in H$, then
{\small \begin{align}
\label{cuspsingularconormal}
N_H(A)(\tau) &= \sum\limits_{\substack{(\alpha,\beta) \in \N_0^n \\ |\alpha| + |\beta| \leq m}}
a_{\alpha,\beta}(x_1,\ldots,x_{k-1},0,y)(x_1^2D_{x_1})^{\alpha_1}\ldots(x_{k-1}^2D_{x_{k-1}})^{\alpha_{k-1}}\tau^{\alpha_k}D_y^{\beta} \\
\intertext{{\normalsize if $p \notin \regbdry\Mbar$, or}}
\label{cuspregularconormal}
N_H(A)(\tau) &= \sum\limits_{\substack{(\alpha,\beta) \in \N_0^n \\ |\alpha| + |\beta| \leq m}}
a_{\alpha,\beta}(x_1,\ldots,x_{k-1},0,y)D_{x_1}^{\alpha_1}(x_2^2D_{x_2})^{\alpha_2}\ldots(x_{k-1}^2D_{x_{k-1}})^{\alpha_{k-1}}\tau^{\alpha_k}
D_y^{\beta}
\end{align}}
if $p \in \regbdry\Mbar$, are local representations of the conormal symbol $N_H(A)(\tau)$
associated with $H$. Observe that $H\cap\regbdry\Mbar$ is the regular part of the boundary
of the manifold with corners $H$.

Everything that we have just said about scalar cusp differential operators holds for operators
acting in sections of vector bundles over $\Mbar$. Let us summarize the above and fix some
notation in the following

\begin{definition}\label{principalsymbolscuspdifferential}
Let $E$ and $F$ be smooth vector bundles over $\Mbar$. The space of cusp differential operators
of order $m \in \N_0$ acting in sections of the bundles $E$ and $F$ will be denoted by
$\Diff_{\cu}^m(\Mbar;E,F)$. By convention, if $\Mbar$ is a closed manifold (i.e. a manifold
with corners of codimension zero), let $\Diff_{\cu}^m(\Mbar;E,F)$ be the space of all
differential operators of order $m$.

Associated with every $A \in \Diff_{\cu}^m(\Mbar;E,F)$ there are the following principal symbols:
\begin{itemize}
\item The \emph{cusp-principal symbol}
$$
\cusym(A) \in C^{\infty}\bigl(\cuT^*\Mbar\setminus 0,
\Hom(\cupi^*E,\cupi^*F)\bigr),
$$
which is a homogeneous function of degree $m \in \N_0$ in the fibres.

Here $\cupi : \cuT^*\Mbar\setminus 0 \to \Mbar$ denotes the canonical projection.
\item The \emph{cusp-principal boundary symbol}
$$
\cusymb(A) \in C^{\infty}\bigl(\cuT^*\regbdry\Mbar \setminus 0,
\Hom(\Scu\otimes\cupi^*E|_{\regbdry\Mbar},\Scu\otimes\cupi^*F|_{\regbdry\Mbar})\bigr),
$$
which is $\kappa$-homogeneous of degree $m \in \N_0$ in the fibres.

Here $\cupi : \cuT^*\regbdry\Mbar \setminus 0 \to \regbdry\Mbar$ is the canonical projection.
\item For each hypersurface $H \subset \singbdry\Mbar$ we have a \emph{conormal symbol} (or normal
operator)
$$
N_H(A)(\tau) \in \Diff_{\cu}^m(H;E|_H,F|_H),
$$
a family of cusp differential operators on $H$ depending on the parameter $\tau \in \R$.
\end{itemize}
\end{definition}

Let $\m$ be a cusp measure on $\Mbar$, i.e. $\x_{\sing}^{(2,\ldots,2)}\m$ is a smooth everywhere
positive density. For any (hermitian) bundle $E \to \Mbar$ let $\L2cu(\Mbar,E)$ be the
$L^2$-space associated with $\m$.

\begin{definition}\label{cuspsobolevspace}
The \emph{cusp Sobolev space} $\Hcu^s(\Mbar,E)$ of sections of $E$ of smoothness $s \in \N_0$
consists of all distributions $u \in \D'(M,E)$ such that $Au \in \L2cu(\Mbar,E)$ for all cusp differential
operators $A$ of order $\leq s$, and let $H_{\cu,0}^s(\Mbar,E)$ be the closure of
$C_0^{\infty}(M,E)$ in $\Hcu^s(\Mbar,E)$.

For $s \in -\N$ we define $\Hcu^s(\Mbar,E)$ as the dual space of $H_{\cu,0}^{-s}(\Mbar,E)$ with
respect to the pairing induced by the $\L2cu$-inner product, and correspondingly let $H_{\cu,0}^{s}(\Mbar,E)$
be the dual space of $\Hcu^{-s}(\Mbar,E)$. The spaces $H_{\cu,0}^s(\Mbar,E)$ and
$\Hcu^s(\Mbar,E)$ for general $s \in \R$ are then defined by interpolation.
\end{definition}

Note that in the case $\regbdry\Mbar = \emptyset$ of standard cusp operators the cusp
Sobolev spaces $H_{\cu,0}^s(\Mbar,E)$ and $\Hcu^s(\Mbar,E)$ coincide, i.e. they differ only near
$\regbdry\Mbar$ in our situation. It is convenient to consider also weighted spaces
$\x_{\sing}^{\alpha}\Hcu^s(\Mbar,E)$ for weights $\alpha = (\alpha_1,\ldots,\alpha_{N-\ell}) \in \R^{N-\ell}$
(recall that $\partial\Mbar$ consists of $\ell$ regular and $N-\ell$ singular
hypersurfaces). By the Sobolev embedding theorem, the space $\dot{C}^{\infty}(\Mbar,E)$ of all $C^{\infty}$-functions
on $\Mbar$ which vanish to infinite order on $\singbdry\Mbar$ equals
\begin{equation}\label{cdotinfty}
\dot{C}^{\infty}(\Mbar,E) = \bigcap\limits_{s \in \R,\; \alpha \in \R^{N-\ell}} \x_{\sing}^{\alpha}\Hcu^s(\Mbar,E),
\end{equation}
and this space is dense in $\x_{\sing}^{\alpha}\Hcu^s(\Mbar,E)$ for all $\alpha \in \R^{N-\ell}$ and $s \in \R$.

For $\alpha,\beta \in \R^{N-\ell}$ with $\alpha_j > \beta_j$ for all $j = 1,\ldots,N-\ell$ and $s > t$ the embedding
\begin{equation}\label{Sobspaceembedding}
\x_{\sing}^{\alpha}H_{\cu(,0)}^s(\Mbar,E) \hookrightarrow \x_{\sing}^{\beta}H_{\cu(,0)}^t(\Mbar,E)
\end{equation}
is compact.

Every cusp operator $A \in \Diff_{\cu}^m(\Mbar;E,F)$ induces continuous operators
$$
A : \x_{\sing}^{\alpha}\Hcu^s(\Mbar,E) \to \x_{\sing}^{\alpha}\Hcu^{s-m}(\Mbar,F)
$$
for all $s \in \R$ and weights $\alpha \in \R^{N-\ell}$. However, when
constructing parametrices and considering therefore pseudodifferential operators, it is
necessary to work with the operator convention $A^+u = r^+Ae^+u$, where $e^+$ denotes the operator
that extends a function $u$ on $M$ by zero to a small neighborhood of
the smooth boundary $\regbdry\Mbar\setminus\singbdry\Mbar$ (observe that there is
a collar neighborhood of the boundary in view of Section \ref{subsec-Manifoldswithcorners}),
then apply an extension of $A$ on this neighborhood to the function $e^+u$,
and finally restrict $Ae^+u$ again via the restriction operator $r^+$ to the interior $M$ of $\Mbar$.
For differential operators we obviously have $A^+ = A$, but formally the operator $e^+$
is well defined only for distributions of Sobolev smoothness $> -\frac{1}{2}$ up to the
smooth part of the boundary, i.e. we have
$$
A^+ : \x_{\sing}^{\alpha}\Hcu^s(\Mbar,E) \to \x_{\sing}^{\alpha}\Hcu^{s-m}(\Mbar,F)
$$
for all $s > -\frac{1}{2}$ and weights $\alpha \in \R^{N-\ell}$.

\subsection{Elliptic boundary problems for cusp differential operators}\label{subsec-Diffboundaryproblems}

Throughout this section let $A \in \Diff_{\cu}^m(\Mbar;E,F)$ be a cusp differential operator.

\begin{definition}\label{intcuspellipt}
$A$ is called \emph{cusp-elliptic} if its cusp-principal symbol $\cusym(A)$ is invertible
on $\cuT^*\Mbar \setminus 0$.
\end{definition}

Let us assume henceforth that $A$ is cusp-elliptic.
An immediate consequence of standard results for ordinary differential equations
is the following

\begin{proposition}\label{cuspbdrysymbsurjective}
The cusp-principal boundary symbol
$$
\cusymb(A) \in C^{\infty}\bigl(\cuT^*\regbdry\Mbar \setminus 0,
\Hom(\Scu\otimes\cupi^*E|_{\regbdry\Mbar},\Scu\otimes\cupi^*F|_{\regbdry\Mbar})\bigr)
$$
is pointwise surjective and has finite dimensional kernel.

We denote by $\K \to \cuT^*\regbdry\Mbar\setminus 0$ the bundle of kernels of $\cusymb(A)$.
\end{proposition}

For any sufficiently smooth section $u$ of a bundle $F$ on $\Mbar\setminus\singbdry\Mbar$
we denote by $\gamma u$ its restriction to the boundary $\regbdry\Mbar \setminus \singbdry\Mbar$,
which gives rise to the restriction operator
\begin{equation}\label{restop}
\gamma : \Hcu^s(\Mbar,F) \to \Hcu^{s-\frac{1}{2}}(\regbdry\Mbar,F|_{\regbdry\Mbar})
\end{equation}
for $s > \frac{1}{2}$. Note that $\Hcu^{s-\frac{1}{2}}(\regbdry\Mbar,F|_{\regbdry\Mbar})$
is just the ordinary cusp Sobolev space of smoothness $s - \frac{1}{2}$
on the manifold with corners $\regbdry\Mbar$ (which coincides with the standard Sobolev
space on all hypersurfaces $H \subset \regbdry\Mbar$ that are smooth).
The operator \eqref{restop} has the following principal symbols:

\begin{definition}\label{principalgamma}
The cusp-principal boundary symbol $\cusymb(\gamma)$ of the restriction operator $\gamma$
is the section ${}^{\cu}\gamma_0 \otimes \id_{\cupi^*F|_{\regbdry\Mbar}}$ of
$\Hom(\Scu\otimes\cupi^*F|_{\regbdry\Mbar},\cupi^*F|_{\regbdry\Mbar})$ on
$\cuT^*\regbdry\Mbar \setminus 0$, where ${}^{\cu}\gamma_0 \in \Hom(\Scu,\C)$ is fibrewise
given by evaluation of a function in $\S(\overline{\R}_+)$ at zero.

Let $H \subset \singbdry\Mbar$ be a singular hypersurface of the boundary which has
nontrivial intersection with $\regbdry\Mbar$. Then the conormal symbol $N_H(\gamma)(\tau)$
of the restriction operator $\gamma$ on the manifold with corners $\Mbar$ is by definition
the constant family $N_H(\gamma)(\tau) \equiv \gamma_H$, $\tau \in \R$, where $\gamma_H$ is the restriction
operator for sections of the bundle $F|_H$ on the manifold with corners $H$ to
$\regbdry H\setminus \singbdry H$. Note that the regular part $\regbdry H$ of the boundary of $H$
is given by $H \cap \regbdry\Mbar$, while $\singbdry H = \singbdry\Mbar \cap \partial H$.
Usually we write just $\gamma$ instead of $\gamma_H$, and the corresponding configuration space
and the vector bundle are always self-understood from the context.
\end{definition}

Now let $B_j \in \Diff_{\cu}^{m_j}(\Mbar;E,F_j)$, $j = 1,\ldots,K$, and let $d = \max\limits_{j=1}^K m_j + 1$.
We then consider the boundary value problem
\begin{equation}\label{diffbvp}
\left.\begin{aligned}
Au &= f \textup{ in } M, \\
Tu &= g \textup{ on } \regbdry\Mbar\setminus\singbdry\Mbar,
\end{aligned}\right\}
\end{equation}
where $T = \bigl(\gamma B_1,\ldots,\gamma B_K\bigr)^{\textup{tr}}$ is the vector of
boundary conditions. Observe that the boundary value problem \eqref{diffbvp} gives rise
to a bounded operator
\begin{equation}\label{diffbvpsobsp}
\A = \binom{A}{T} : \Hcu^s(\Mbar,E) \to \begin{array}{c} \Hcu^{s-m}(\Mbar,F) \\ \oplus \\
\bigoplus\limits_{j=1}^K \Hcu^{s-m_j-\frac{1}{2}}(\regbdry\Mbar,F_j|_{\regbdry\Mbar}) \end{array}
\end{equation}
for $s > d-\frac{1}{2}$.

\begin{definition}\label{ellipticitybvp}
Let $A$ be cusp-elliptic. We call the boundary value problem $\A = \binom{A}{T}$ \emph{elliptic},
if the following conditions are fulfilled:
\begin{enumerate}[i)]
\item The mapping
$$
\cusymb(T) = \begin{pmatrix} \cusymb(\gamma)\cusymb(B_1) \\ \vdots \\ \cusymb(\gamma)\cusymb(B_K) \end{pmatrix} :
\K \to \bigoplus\limits_{j=1}^K \cupi^*F_j|_{\regbdry\Mbar}
$$
is a vector bundle isomorphism on $\cuT^*\regbdry\Mbar \setminus 0$. Recall that $\K$ is
the bundle of kernels of $\cusymb(A)$.

This condition is equivalent to the invertibility of the \emph{cusp-principal boundary symbol}
$$
\cusymb(\A) = \binom{\cusymb(A)}{\cusymb(T)} : \Scu\otimes\cupi^*E|_{\regbdry\Mbar} \to \begin{array}{c} \Scu\otimes\cupi^*F|_{\regbdry\Mbar} \\ \oplus \\
\bigoplus\limits_{j=1}^K \cupi^*F_j|_{\regbdry\Mbar} \end{array}
$$
of the boundary value problem $\A$ on $\cuT^*\regbdry\Mbar \setminus 0$, and it is the
appropriate version of the \emph{Shapiro-Lopatinsky condition} in our context of boundary
problems for cusp operators.
\item For each singular boundary hypersurface $H \subset \singbdry\Mbar$ with
$H \cap \regbdry\Mbar = \emptyset$ the conormal symbol
$$
N_H(\A)(\tau) \equiv N_H(A)(\tau) : \Hcu^{s}(H,E|_H) \to \Hcu^{s-m}(H,F|_H)
$$
is invertible for all $\tau \in \R$ and some (all) $s \in \R$.
\item For each singular boundary hypersurface $H \subset \singbdry\Mbar$ with
$H \cap \regbdry\Mbar \neq \emptyset$ the \emph{conormal symbol}
$$
N_H(\A)(\tau) = \binom{N_H(A)(\tau)}{N_H(T)(\tau)} : \Hcu^s(H,E|_H) \to \begin{array}{c}
\Hcu^{s-m}(H,F|_H) \\ \oplus \\ \bigoplus\limits_{j=1}^K \Hcu^{s-m_j-\frac{1}{2}}(\regbdry H,F_j|_{\regbdry H})
\end{array}
$$
of the boundary value problem $\A$ is invertible for all $\tau \in \R$ and
some (all) $s > \max\{m,d\}-\frac{1}{2}$. Here we write analogously to i)
$$
N_H(T)(\tau) = \begin{pmatrix} N_H(\gamma)(\tau) N_H(B_1)(\tau) \\ \vdots \\ N_H(\gamma)(\tau) N_H(B_K)(\tau) \end{pmatrix}.
$$
Observe that $N_H(\A)(\tau)$ is a family of boundary value problems on the manifold with
corners $H$.
\end{enumerate}
\end{definition}

\begin{remark}\label{remconormalsymb}
Assume that $A$ is cusp-elliptic, and the boundary value problem $\A = \binom{A}{T}$ satisfies only condition i) in
Definition \ref{ellipticitybvp}. Then, in view of Theorem \ref{largeparaminverse}, the
conormal symbols $N_H(\A)(\tau)$ in ii) and iii) of Definition \ref{ellipticitybvp} are
automatically invertible for $|\tau| > 0$ sufficiently large, and the inverses $N_H(\A)(\tau)^{-1}$
are represented as families of cusp pseudodifferential operators resp. boundary value problems
depending on the parameter $\tau \in \R$. Therefore, the conditions ii) and iii) in
Definition \ref{ellipticitybvp} are in a sense subordinate to the invertibility of the
cusp-principal symbol and the cusp-principal boundary symbol, but these requirements
nevertheless are essential for the validity of Theorem \ref{diffbvpFredholm} below.
\end{remark}

\begin{theorem}\label{diffbvpFredholm}
Assume that the boundary value problem $\A = \binom{A}{T}$ is elliptic in the sense of
Definition \ref{ellipticitybvp}. Then, for $s > \max\{m,d\}-\frac{1}{2}$ and all $\alpha \in \R^{N-\ell}$, the operator
\begin{equation}\label{bvpsob}
\A = \binom{A}{T} : \x_{\sing}^{\alpha}\Hcu^s(\Mbar,E) \to \begin{array}{c} \x_{\sing}^{\alpha}\Hcu^{s-m}(\Mbar,F) \\ \oplus \\
\bigoplus\limits_{j=1}^K \x_{\sing}^{\alpha}\Hcu^{s-m_j-\frac{1}{2}}(\regbdry\Mbar,F_j|_{\regbdry\Mbar}) \end{array}
\end{equation}
is a Fredholm operator. Recall that $d = \max\limits_{j=1}^K m_j + 1$, where $m_j$ is the order of
the boundary condition $B_j$.

The kernel $N(\A)$ of \eqref{bvpsob} is a subspace of $\dot{C}^{\infty}(\Mbar,E)$, the space of
all $C^{\infty}$-functions on $\Mbar$ which vanish to infinite order on $\singbdry\Mbar$, and
therefore does not depend on $s > \max\{m,d\}-\frac{1}{2}$ and $\alpha \in \R^{N-\ell}$.
More generally, if
$$
(f,g) \in \x_{\sing}^{\alpha}\Hcu^{s-m}(\Mbar,F) \oplus \bigoplus\limits_{j=1}^K \x_{\sing}^{\alpha}\Hcu^{s-m_j-\frac{1}{2}}(\regbdry\Mbar,F_j|_{\regbdry\Mbar})
$$
for $s > \max\{m,d\}-\frac{1}{2}$ and $\alpha \in \R^{N-\ell}$, and if the function $u \in \x_{\sing}^{\beta}\Hcu^t(\Mbar,E)$ is
a solution of $Au = f$ and $Tu = g$ for some $t > \max\{m,d\}-\frac{1}{2}$ and some $\beta \in \R^{N-\ell}$, then $u \in \x_{\sing}^{\alpha}\Hcu^s(\Mbar,E)$.

There exists a parametrix
$$
{\mathcal P} = \begin{pmatrix} P^+ + G & K_1 & \cdots & K_K \end{pmatrix}
$$
of $\A$ in the cusp calculus of pseudodifferential boundary value problems that is defined
in Section \ref{sec-Pseudocalculus}, where $P^+$ has order $-m$, $G$ is a singular Green operator
of order $-m$ and type $(d-m)_+ = \max\{0,d-m\}$, and the $K_j$ are generalized singular Green
operators of potential type of order $-m_j-\frac{1}{2}$ (see Theorem \ref{Parametrix}).
\end{theorem}
\begin{proof}
We employ the pseudodifferential calculus from Section \ref{sec-Pseudocalculus}.
According to Lemma \ref{Orderreductions} there exist order reductions (i.e. elliptic
invertible operators)
$$
R_j : \x_{\sing}^{\alpha}\Hcu^s(\regbdry\Mbar,F_j|_{\regbdry\Mbar}) \to \x_{\sing}^{\alpha}\Hcu^{s-\mu_j}(\regbdry\Mbar,F_j|_{\regbdry\Mbar})
$$
in the class $\Psicu^{\mu_j}(\regbdry\Mbar)$ of
cusp operators on $\regbdry\Mbar$, where $\mu_j = m-m_j-\frac{1}{2}$, $j = 1,\ldots,K$.
Let
$$
R = \begin{pmatrix} R_1 & 0 & 0 \\ \vdots & \ddots & \vdots \\ 0 & 0 & R_K \end{pmatrix},
$$
and consider for $s > \max\{m,d\}-\frac{1}{2}$ and $\alpha \in \R^{N-\ell}$ the operator
$$
{\mathcal B} = \begin{pmatrix} 1 & 0 \\ 0 & R \end{pmatrix}\binom{A}{T} = \binom{A}{RT} :
\x_{\sing}^{\alpha}\Hcu^s(\Mbar,E) \to
\begin{array}{c} \x_{\sing}^{\alpha}\Hcu^{s-m}(\Mbar,F) \\ \oplus \\
\x_{\sing}^{\alpha}\Hcu^{s-m}(\regbdry\Mbar,J_+) \end{array},
$$
where $J_+ = \bigoplus\limits_{j=1}^KF_j|_{\regbdry\Mbar}$.
Hence, at the expense of changing the differential boundary condition $T$ to the pseudodifferential
boundary condition $RT$, we have obtained an elliptic pseudodifferential boundary value problem
${\mathcal B}$ in the class $\Psicu^{m,d}(\Mbar)$ defined in Section \ref{sec-Pseudocalculus} --- that
${\mathcal B}$ is indeed elliptic (Definition \ref{pseudoelliptic}) follows from the
multiplicativity of the principal symbols under composition and our ellipticity assumption
about $\A$ according to Definition \ref{ellipticitybvp}.

By Theorem \ref{Parametrix} there exists a parametrix ${\mathcal P}' \in \Psicu^{-m,(d-m)_+}(\Mbar)$
of ${\mathcal B}$ up to residual Green operators in the class $\Psicu^{-\infty,*}_{\infty}(\Mbar)$, and consequently
$$
{\mathcal P} = {\mathcal P}'\begin{pmatrix} 1 & 0 \\ 0 & R \end{pmatrix}
$$
is a parametrix of $\A$ as desired.
The Fredholmness and elliptic regularity follow from the mapping properties of ${\mathcal P}$ and
of the remainders, see also Corollary \ref{elliptFredh} in Section \ref{sec-Pseudocalculus}.
\end{proof}

In the proof of Theorem \ref{diffbvpFredholm} we have used the following lemma, which is a
standard result in the theory of pseudodifferential operators, and it also
holds within the ordinary cusp algebra.

\begin{lemma}\label{Orderreductions}
For every vector bundle $F$ and every $\mu \in \R$ there exists a reduction of orders, i.e.
an elliptic invertible operator
$$
R : \x_{\sing}^{\alpha}\Hcu^s(\regbdry\Mbar,F) \to \x_{\sing}^{\alpha}\Hcu^{s-\mu}(\regbdry\Mbar,F)
$$
in the class of cusp operators $\Psicu^{\mu}(\regbdry\Mbar)$ of order $\mu$, and the inverse
$R^{-1}$ belongs to the class $\Psicu^{-\mu}(\regbdry\Mbar)$.
\end{lemma}
\begin{proof}
Let $R(\tau) \in \Psicu^{\mu}(\regbdry\Mbar,\R)$ be any parameter-dependent cusp-elliptic operator
in the class of cusp operators of order $\mu$ on $\regbdry\Mbar$ which
depend on the parameter $\tau \in \R$ (strong polyhomogeneity, i.e. locally modelled on
classical symbols).
Such an operator always exists because there exist elliptic (i.e. invertible) parameter-dependent
cusp-principal symbols on $\bigl(\cuT^*\regbdry\Mbar \times \R\bigr) \setminus 0$ (just choose
a metric), and any cusp-quantization of such a symbol yields an operator $R(\tau)$ as stated.

Thus, via constructing a parametrix in the class of parameter-dependent cusp operators
(see also Theorem \ref{largeparaminverse}), we obtain that $R(\tau)$ is invertible for
$|\tau| > 0$ sufficiently large, and the inverse $R(\tau)^{-1} \in \Psicu^{-\mu}(\regbdry\Mbar,\R)$
(more precisely, there is a parameter-dependent parametrix of $R(\tau)$ in
$\Psicu^{-\mu}(\regbdry\Mbar,\R)$ which equals the inverse for large $\tau$).

Hence we can choose $R = R(\tau_0)$ for some $\tau_0 \in \R$ with $|\tau_0| > 0$ sufficiently large.
\end{proof}


\section{The cusp calculus of pseudodifferential boundary value problems}\label{sec-Pseudocalculus}

The aim of this section is to set up a Boutet de Monvel's calculus of cusp pseudodifferential
boundary value problems on the manifold with corners $\Mbar$. We assume that the reader is
familiar with the classical Boutet de Monvel algebra on an arbitrary smooth
manifold with boundary. Otherwise the entries \cite{GrubbBuch, GruSe95, KaSchu03, KraiParabBVP,
SchroSchu, SzWiley98} in the list of references will provide useful introductory information.

In order to understand the conormal symbolic structure associated with singular hypersurfaces
$H \subset \singbdry\Mbar$ --- see \eqref{cuspsingularconormal}, \eqref{cuspregularconormal}, and
Definitions \ref{principalgamma} and \ref{ellipticitybvp} --- it is necessary to admit from the
very beginning that the operators depend on a parameter $\lambda \in \Lambda$, where
$\Lambda \subset \R^q$ is the closure of any open conical subset of $\R^q$, or $\Lambda = \{0\}$
(the case without parameters). This is of course of independent interest also, e.g., for the analysis
of resolvents of elliptic boundary value problems and the heat equation.

Initially, the operator families $A(\lambda)$ in the cusp calculus
\begin{equation}\label{psicuast}
\Psicu^{\ast,\ast}(\Mbar,\Lambda) = \bigcup\limits_{\substack{\mu \in \Z \\ d \in \N_0}}\Psicu^{\mu,d}(\Mbar,\Lambda)
\end{equation}
are considered in the spaces
\begin{equation}\label{initialcusp}
A(\lambda) : \begin{array}{c} \dot{C}^{\infty}(\Mbar,E) \\ \oplus \\ \dot{C}^{\infty}(\regbdry\Mbar,J_-) \end{array}
\to \begin{array}{c} \dot{C}^{\infty}(\Mbar,F) \\ \oplus \\ \dot{C}^{\infty}(\regbdry\Mbar,J_+) \end{array}.
\end{equation}
Here $E$ and $F$ are smooth vector bundles over $\Mbar$, and $J_{\pm}$ are smooth vector bundles
over $\regbdry\Mbar$. Note that the vector bundles $J_+$ or $J_-$ are admitted to be zero which
happens to be the case, in particular, for differential boundary value problems and their
parametrices, see Section \ref{subsec-Diffboundaryproblems}. Recall that $\dot{C}^{\infty}$ always
denotes the space of smooth functions on a manifold with corners that vanish to infinite order
on the singular part of the boundary.

The cusp algebra \eqref{psicuast} ist filtered by the pseudodifferential order
$\mu \in \Z$, while $d \in \N_0$ is the type of the pseudodifferential boundary value problem in Boutet de Monvel's algebra.
There is a second filtration by weights, namely we consider for $\alpha = (\alpha_1,\ldots,\alpha_{N-\ell})
\in \N_0^{N-\ell}$ the ideals $\x_{\sing}^{\alpha}\Psicu^{\ast,\ast}(\Mbar,\Lambda)$
which encode the order of vanishing separately on each hypersurface $H \subset \singbdry\Mbar$.
When considering cusp operators in the spaces \eqref{initialcusp}, we could even admit general
real weights $\alpha \in \R^{N-\ell}$, but we would leave the class \eqref{psicuast} then. Recall that,
by general convention in pseudodifferential boundary value problems,
the matrix multiplication operator by a smooth function $\varphi \in C^{\infty}(\overline{N})$ is given by
\begin{equation}\label{Multop}
\varphi \equiv M_{\varphi} = \begin{pmatrix} \varphi & 0 \\ 0 & \varphi|_{\partial\overline{N}} \end{pmatrix} :
\begin{array}{c} C^{\infty}(\overline{N},E) \\ \oplus \\ C^{\infty}(\partial\overline{N},J_-) \end{array}
\to \begin{array}{c} C^{\infty}(\overline{N},E) \\ \oplus \\ C^{\infty}(\partial\overline{N},J_-) \end{array},
\end{equation}
where $\overline{N} = \Mbar \setminus \singbdry\Mbar$, and $\partial\overline{N} = \partial\Mbar \setminus
\singbdry\Mbar$, form the regular part of $\Mbar$ and $\partial\Mbar$. In particular, multiplication of an
operator $A(\lambda) \in \Psicu^{\ast,\ast}(\Mbar,\Lambda)$ with the function $\x_{\sing}^{\alpha}$
makes sense, and in this way the ideals $\x_{\sing}^{\alpha}\Psicu^{\ast,\ast}(\Mbar,\Lambda)$ are
defined.

Our construction of $\Psicu^{\ast,\ast}(\Mbar,\Lambda)$ is performed in symbolical terms, and we proceed
by induction on the codimension of the manifold with corners. In the coordinates $x = \x_i$ transversal
to hypersurfaces $H \subset \singbdry\Mbar$, the quantization makes use of the
cusp transform
\begin{align}
\F_{\cu} &: C_0^{\infty}(\R_+) \to \S(\R), \\
\label{cusptrafo}
\bigl(\F_{\cu}u\bigr)(\xi) &= \int\limits_0^{\infty}e^{i\xi/x}u(x)\,\frac{dx}{x^2}, \\
\intertext{and the inverse cusp transform}
\F_{\cu}^{-1} &: C_0^{\infty}(\R) \to C^{\infty}(\R_+), \\
\bigl(\F_{\cu}^{-1}u\bigr)(x) &= \frac{1}{2\pi}\int\limits_{\R}e^{-i\xi/x}u(\xi)\,d\xi.
\end{align}
Observe that $\F_{\cu}\bigl(x^2D_xu\bigr) = \xi \F_{\cu}u$, and a change of variables
$t = -\frac{1}{x}$ in \eqref{cusptrafo} reveals 
$\bigl(\F_{\cu}u\bigr)(\tau) = \int\limits_{-\infty}^{0}e^{-it\tau}u(-1/t)\,dt$, i.e.
the cusp transform of $u$ equals the Fourier transform of the function $u(-1/t) \in C_0^{\infty}(-\infty,0)$.
These properties suggest to make use of the cusp transform in our situation.

\medskip

\noindent
We proceed as follows:
\begin{itemize}
\item In Section \ref{subsec-ResidualGreen} we define the class of regularizing residual
Green operators in the cusp algebra, the smallest ideal in the calculus.
\item By induction we assume that we already know the cusp calculus on manifolds with
corners of codimension $\leq \codim\Mbar - 1$. Making use of this for the singular hypersurfaces
$H \subset \singbdry\Mbar$, we construct in Section \ref{subsec-SingularCollar} the class of
cusp operators in every collar neighborhood $[0,\eps)\times H$.
\item Finally, in Section \ref{subsec-Fullalgebra}, we construct the cusp calculus
on $\Mbar$, and by induction we then have the cusp algebra of pseudodifferential boundary value
problems on all manifolds with corners.
\end{itemize}
In order to make sense of this inductive process, note that in the case of a manifold $\Mbar$
with corners of codimension zero, i.e. a closed compact manifold without boundary, we simply let
$\Psicu^{\mu,d}(\Mbar,\Lambda) \equiv L_{\cl}^{\mu}(\Mbar,\Lambda)$ be the standard class of classical
(parameter-dependent) pseudodifferential operators on $\Mbar$ which are locally modelled on
symbols $a(z,\zeta,\lambda)$ such that
$|\partial_z^{\alpha}\partial_{(\zeta,\lambda)}^{\beta}a| = O(|\zeta,\lambda|^{\mu-|\beta|})$
as $|\zeta,\lambda| \to \infty$, uniformly for $z$ in compact
sets, and such that $a$ admits an asymptotic expansion $a \sim \sum\limits_{j=0}^{\infty}a_{(\mu-j)}$
into homogeneous components $a_{(\mu-j)}(z,\zeta,\lambda)$ (with respect to $(\zeta,\lambda) \neq 0$) of degree $\mu-j$.

\subsection{The class $\Psicu_{\infty}^{-\infty,d}(\Mbar,\Lambda)$ of residual Green operators}\label{subsec-ResidualGreen}

The union of all collar neighborhoods of hypersurfaces $H \subset \regbdry\Mbar$ induces, by
assumption about $\regbdry\Mbar$, one collar neighborhood
$[0,\eps)\times\regbdry\Mbar \cong \Mbar$ of $\regbdry\Mbar$ in $\Mbar$.

Let $\partial \in \Diff_{\cu}^1(\Mbar)$ be a cusp vector field on $\Mbar$
supported within $[0,\eps/2)\times\regbdry\Mbar$,
which coincides near $\regbdry\Mbar$ with $\partial_x$, where $x$ is the coordinate in $[0,\eps)$,
and consider the operator
$$
\partial_+ = r^+\partial e^+ : \Hcu^s(\Mbar) \to {\mathscr D}'(M), \quad s > -\frac{1}{2},
$$
defined by extending a distribution $u \in \Hcu^s(\Mbar)$ by zero to a small cylinder $(-\delta,\eps/2)\times\regbdry\Mbar$
around $\regbdry\Mbar$, differentiating the resulting distribution with respect to
$x \in (-\delta,\eps/2)$ (i.e. applying the canonical extension of $\partial$ to it), and
restricting it again to the interior $M$ of $\Mbar$ (see also the remarks at the end of
Section \ref{subsec-cuspdifferentialoperators}). By identifying a vector bundle $E$ with
the pull-back of its restriction to the boundary in the collar neighborhood
$[0,\eps)\times\regbdry\Mbar$ this operation extends in an obvious way to sections of $E$ over
$\Mbar$, and the compositions define continuous operators
$$
\partial_+^j : \Hcu^s(\Mbar,E) \to \Hcu^{s-j}(\Mbar,E)
$$
for $j \in \N_0$ provided that $s > j - \frac{1}{2}$.
Counting these derivatives constitute the type of an operator in Boutet de Monvel's
calculus (which at the end does not depend on all the choices involved here).

Let us first define the class $\Psicu_{\infty}^{-\infty,0}(\Mbar,\Lambda)$ of residual Green
operators of type zero. An operator (or a family of operators) $G(\lambda)$ belongs to this class
if and only if, for all $k \in \N_0$ and multiindices $|\alpha| \leq k$,
$$
G(\lambda) : \begin{array}{c} \L2cu(\Mbar,E) \\ \oplus \\ \L2cu(\regbdry\Mbar,J_-) \end{array} \to \begin{array}{c} \x_{\textup{sing}}^{\alpha}\Hcu^k(\Mbar,F) \\ \oplus \\ \x_{\textup{sing}}^{\alpha}\Hcu^k(\regbdry\Mbar,J_+) \end{array}
$$
is continuous and rapidly decreasing with all derivatives as $|\lambda| \to \infty$
in $\Lambda$ (where the latter condition is void if $\Lambda$ is just a point), and
correspondingly so for the formal adjoint
$$
G(\lambda)^* : \begin{array}{c} \L2cu(\Mbar,F) \\ \oplus \\ \L2cu(\regbdry\Mbar,J_+) \end{array} \to \begin{array}{c} \x_{\textup{sing}}^{\alpha}\Hcu^k(\Mbar,E) \\ \oplus \\ \x_{\textup{sing}}^{\alpha}\Hcu^k(\regbdry\Mbar,J_-) \end{array}
$$
with respect to the pairing induced by the $\bigl(\L2cu \oplus \L2cu\bigr)$-inner products.
Recall that the cusp Sobolev spaces $\Hcu^s$ coincide
with the standard Sobolev spaces on all smooth hypersurfaces of $\regbdry\Mbar$.
Consequently, the class $\Psicu_{\infty}^{-\infty,0}(\Mbar,\Lambda)$ consists of all operators
with smooth kernels that vanish to infinite order at $\singbdry\Mbar$ and depend rapidly
decreasing with all derivatives on the parameter $\lambda \in \Lambda$.

In general an operator $G(\lambda) \in \Psicu_{\infty}^{-\infty,d}(\Mbar,\Lambda)$ is a residual
Green operator of type $d \in \N_0$, if it can be written as a sum
$$
G(\lambda) = \sum\limits_{j=0}^d G_j(\lambda)\begin{pmatrix} \partial_+ & 0 \\ 0 & 0 \end{pmatrix}^j :
\begin{array}{c} \Hcu^s(\Mbar,E) \\ \oplus \\ \Hcu^s(\regbdry\Mbar,J_-) \end{array} \to
\begin{array}{c} \Hcu^{s'}(\Mbar,F) \\ \oplus \\ \Hcu^{s'}(\regbdry\Mbar,J_+) \end{array}
$$
for $s > d - \frac{1}{2}$, where the $G_j(\lambda) \in \Psicu_{\infty}^{-\infty,0}(\Mbar,\Lambda)$ are residual Green operators of type zero.
Note that the class $\Psicu_{\infty}^{-\infty,d}(\Mbar,\Lambda)$ carries a
natural Fr{\'e}chet topology as a nondirect sum of the Fr{\'e}chet spaces $\Psicu_{\infty}^{-\infty,0}(\Mbar,\Lambda)$,
which themselves carry the topology induced by the kernels or, equivalently, by the defining mapping properties
stated above with rapidly decreasing dependence on $\lambda \in \Lambda$ (with all derivatives).

\subsection{The cusp calculus near a singular boundary hypersurface}\label{subsec-SingularCollar}

By induction we now assume that we know the cusp calculus of pseudodifferential boundary
value problems on manifolds with corners of codimension $\leq \codim\Mbar - 1$. Let
$H \subset \singbdry\Mbar$ be a singular boundary hypersurface of $\Mbar$, and let
$[0,\eps)\times H \cong \Mbar$ be the collar neighborhood associated with $H$ and the defining
function $\x$ for $H$.

$H$ itself is a manifold with corners of codimension $\leq \codim\Mbar - 1$, and
$\regbdry H = H \cap \regbdry\Mbar$ is the regular part of the boundary $\partial H$ of $H$,
while $\singbdry H = \partial H \cap \singbdry\Mbar$ is the singular part.

Our setting makes it necessary to consider two cases: The first and essential case
$\regbdry H \neq \emptyset$, and the second case $\partial H = \singbdry H$ of
ordinary cusp operators. We focus in the sequel on $\regbdry H \neq \emptyset$, the case
$\regbdry H = \emptyset$ is simpler (just ignore all boundary related constructions below).
A not symbolical, but kernel-oriented definition of the ordinary cusp algebra
(the case $\regbdry\Mbar = \emptyset$) can be found in \cite{LauterMoroianu}.

\begin{definition}\label{cuspcollardef}
An operator family
$$
A(\lambda) : \dot{C}^{\infty}_0\biggl([0,\eps),\begin{array}{c} \dot{C}^{\infty}(H,E) \\ \oplus \\ \dot{C}^{\infty}(\regbdry H,J_-) \end{array}\biggr) \to
\dot{C}^{\infty}\biggl([0,\eps),\begin{array}{c} \dot{C}^{\infty}(H,F) \\ \oplus \\ \dot{C}^{\infty}(\regbdry H,J_+) \end{array}\biggr)
$$
belongs to the space $\Psicu^{\mu,d}([0,\eps)\times H,\Lambda)$
of parameter-dependent cusp operators of order $\mu \in \Z$ and type $d \in \N_0$, if it is of
the form
\begin{equation}\label{cuspquantization}
\bigl(A(\lambda)u\bigr)(x) = \frac{1}{2\pi}\int\limits_{\R}\int\limits_{0}^{\eps}
e^{i(1/y-1/x)\xi}a(x,\xi,\lambda)u(y)\,\frac{dy}{y^2}\,d\xi + \bigl(C(\lambda)u\bigr)(x),
\end{equation}
where
\begin{gather}
\label{Cdef}\bigl(C(\lambda)u\bigr)(x) = \int\limits_0^{\eps}c(x,y,\lambda)u(y)\,\frac{dy}{y^2}, \\
\label{opkernel}c(x,y,\lambda) \in \dot{C}^{\infty}\bigl([0,\eps)_x\times[0,\eps)_y,\Psicu^{-\infty,d}(H,\Lambda)\bigr), \\
\intertext{and}
\label{cuspsymbol}a(x,\xi,\lambda) \in C^{\infty}\bigl([0,\eps)_x,\Psicu^{\mu,d}(H,\R_{\xi}\times\Lambda)\bigr).
\end{gather}
Note that we know by induction the class $\Psicu^{\mu,d}(H,\R\times\Lambda)$ of
cusp pseudodifferential boundary value problems on $H$ depending on the parameters $(\xi,\lambda) \in \R\times\Lambda$,
and this space is endowed with a natural Fr{\'e}chet topology (this topology is also known
by induction). Consequently, the spaces of operator-valued kernels \eqref{opkernel} and
symbols \eqref{cuspsymbol} are well defined and carry themselves natural Fr{\'e}chet topologies, and
by \eqref{cuspquantization} we thus also have a natural Fr{\'e}chet topology on
$\Psicu^{\mu,d}([0,\eps)\times H,\Lambda)$ as a nondirect sum of Fr{\'e}chet spaces.
\end{definition}

For $u \in C^{\infty}_0\biggl((0,\eps),\begin{array}{c} \dot{C}^{\infty}(H,E) \\ \oplus \\ \dot{C}^{\infty}(\regbdry H,J_-) \end{array}\biggr)$
a change of variables in \eqref{cuspquantization} reveals

\begin{align*}
\bigl(A(\lambda)u\bigr)(x) &= \frac{1}{2\pi}\int\limits_{\R}\int\limits_{0}^{\eps}
e^{i(x-y)\xi}a(x,xy\xi,\lambda)\frac{x}{y}u(y)\,dy\,d\xi + \bigl(C(\lambda)u\bigr)(x) \\
&= \frac{1}{2\pi}\int\limits_{\R}\int\limits_{0}^{\eps}
e^{i(x-y)\xi}a(x,x^2\xi,\lambda)u(y)\,dy\,d\xi + \bigl(B(\lambda)u\bigr)(x),
\end{align*}
where $B(\lambda)$ is a parameter-dependent operator in Boutet de Monvel's calculus
of order $\mu-1$ and type $d$ on the regular part of $(0,\eps)\times H$. From this identity
and the induction hypothesis, we see that in coordinates $[0,\eps)^k\times\Omega$
on $[0,\eps)\times H$, where all coordinates $(x_1,\ldots,x_k) \in [0,\eps)^k$ are associated
with singular hypersurfaces at $x_i = 0$, every operator $A(\lambda) \in
\Psicu^{\mu,d}([0,\eps)\times H,\Lambda)$ has a homogeneous principal symbol
of the form $\sigma(x_1,\ldots,x_k,y,x_1^2\xi_1,\ldots,x_k^2\xi_k,\eta,\lambda)$ with a
homogeneous function
$$
\sigma(x_1,\ldots,x_k,y,\xi_1,\ldots,\xi_k,\eta,\lambda)
$$
in the variables $(\xi_1,\ldots,\xi_k,\eta,\lambda) \neq 0$ of degree $\mu$,
and correspondingly so for the principal boundary symbol (which is $\kappa$-homogeneous of degree $\mu$).
Note that $A(\lambda)$ is, in particular, an operator in Boutet de Monvel's calculus on the regular part of $(0,\eps)\times H$.
Consequently, the principal symbol and the principal boundary symbol of every operator matrix
$A(\lambda) \in \Psicu^{\mu,d}([0,\eps)\times H,\Lambda)$ extend to well
defined sections on the cusp cotangent bundles
$\bigl(\cuT^*\bigl([0,\eps)\times H\bigr) \times \Lambda\bigr)\setminus 0$
and $\bigl(\cuT^*\bigl([0,\eps)\times \regbdry H\bigr) \times \Lambda\bigr) \setminus 0$, respectively.
The so obtained symbols are the cusp-principal symbol $\cusym(A)$ and the cusp-principal boundary symbol
$\cusymb(A)$ of $A(\lambda)$, see also Definition \ref{principalsymbolcusppseudo} in the context of the
full cusp calculus.

\begin{proposition}\label{cusplocalsplitexact}
The cusp-principal symbol sequence
$$
\begin{CD}
0 @>>> \Psicu^{\mu-1,d} @>>> \Psicu^{\mu,d} @>>(\cusym,\cusymb)> {}^{\cu}\Sigma @>>> 0
\end{CD}
$$
for (parameter-dependent) operators on $[0,\eps)\times H$ is topologically split exact.

${}^{\cu}\Sigma$ -- the space of principal symbols -- consists of tuples of homogeneous (resp. $\kappa$-homogeneous)
sections on the cusp cotangent bundles that satisfy a canonical compatibility condition.
\end{proposition}
\begin{proof}
By induction the cusp-principal symbol sequence for operators on $H$ depending on the parameters
$(\xi,\lambda) \in \R\times\Lambda$ is topologically split exact. Let
$$
{}^{\cu}\textup{op} : {}^{\cu}\Sigma(H,\R\times\Lambda) \to \Psicu^{\mu,d}(H,\R\times\Lambda)
$$
be any quantization map, i.e. any continuous right inverse of the cusp-principal symbol mapping.
Every element $\sigma \in {}^{\cu}\Sigma([0,\eps)\times H,\Lambda)$ can be represented in a unique way as
$\sigma(x,x^2\xi,\lambda)$ with $\sigma(x,\xi,\lambda) \in C^{\infty}\bigl([0,\eps)_x,{}^{\cu}\Sigma(H,\R_{\xi}\times\Lambda)\bigr)$,
and so
$$
\sigma \mapsto a_{\sigma}(x,\xi,\lambda) = {}^{\cu}\textup{op}\,\sigma(x,\cdot,\cdot) \in C^{\infty}\bigl([0,\eps),\Psicu^{\mu,d}(H,\R\times\Lambda)\bigr)
$$
is a continuous linear mapping which associates with a tuple $\sigma$ of cusp-principal symbols
an operator-valued symbol of the form \eqref{cuspsymbol}. We now define the operator
$A_{\sigma}(\lambda) \in \Psicu^{\mu,d}([0,\eps)\times H,\Lambda)$ according to
\eqref{cuspquantization} with the operator-valued symbol $a(x,\xi,\lambda) = a_{\sigma}(x,\xi,\lambda)$ and
$C(\lambda) \equiv 0$. This gives rise to a quantization mapping
$$
{}^{\cu}\Sigma([0,\eps)\times H,\Lambda) \ni \sigma \mapsto A_{\sigma}(\lambda) \in
\Psicu^{\mu,d}([0,\eps)\times H,\Lambda)
$$
as desired.
\end{proof}

\begin{notation}
For functions $\varphi$ and $\psi$ we write $\varphi \prec \psi$ if $\psi \equiv 1$ in a neighborhood
of the support of $\varphi$.
\end{notation}

\begin{lemma}\label{localpseudolocal}
Let $\omega \prec \tilde{\omega}$ be cut-off functions near zero, i.e. $\omega,\tilde{\omega} \in C_0^{\infty}([0,\eps))$
with $\omega,\tilde{\omega} \equiv 1$ near $x = 0$. Furthermore, let $A(\lambda) \in \Psicu^{\mu,d}([0,\eps)\times H,\Lambda)$.

Then the operators $\omega A(\lambda) (1-\tilde{\omega})$ and $(1-\tilde{\omega}) A(\lambda) \omega$
are of the form \eqref{Cdef} with kernels \eqref{opkernel}.
\end{lemma}
\begin{proof}
Let $\tilde{A}(\lambda)$ be any of the operators $\omega A(\lambda) (1-\tilde{\omega})$ or $(1-\tilde{\omega}) A(\lambda) \omega$.
Obviously, $\tilde{A}(\lambda)$ is of the form \eqref{Cdef} if and only if $\hat{\omega} \tilde{A}(\lambda) \hat{\omega}$
is such for all cut-off functions $\hat{\omega} \in C_0^{\infty}([0,\eps))$ near zero. However, the operator
$\hat{\omega} \tilde{A}(\lambda) \hat{\omega}$ is of the form $\varphi A(\lambda) \psi$ for functions
$\varphi, \psi \in C_0^{\infty}([0,\eps))$ with $\supp \varphi \cap \supp \psi = \emptyset$. In coordinates
$t = -\frac{1}{x}$, $\varphi A(\lambda) \psi$ thus takes the form
\begin{equation}\label{transformedoperator}
\frac{1}{2\pi}\iint e^{i(t-t')\tau}b(t,t',\tau,\lambda)u(t')\,dt'\,d\tau
\end{equation}
with a symbol $b(t,t',\tau,\lambda) \in S^{0,0}(\R_t\times\R_{t'},\Psicu^{\mu,d}(H,\R_{\tau}\times\Lambda))$
that vanishes for $|t-t'| < \delta$. Consequently, by oscillatory integral techniques for pseudodifferential operators
with global symbols (see Remark \ref{transformation} below), this operator is an integral operator with operator-valued kernel
in $\S(\R_t\times\R_{t'},\Psicu^{-\infty,d}(H,\Lambda))$ in the $t$-coordinates. Transforming back to the $x$-coordinates,
this is just what is claimed.
\end{proof}

\begin{remark}\label{transformation}
In the proof of Lemma \ref{localpseudolocal} we employed an observation which is very helpful also further below
for setting up the properties of cusp pseudodifferential operators:

Let $\varphi,\psi \in C_0^{\infty}([0,\eps))$ and $A(\lambda) \in \Psicu^{\mu,d}([0,\eps)\times H,\Lambda)$, and
consider the operator
$$
\varphi A(\lambda) \psi : \dot{C}^{\infty}_0\biggl([0,\eps),\begin{array}{c} \dot{C}^{\infty}(H,E) \\ \oplus \\ \dot{C}^{\infty}(\regbdry H,J_-) \end{array}\biggr) \to
\dot{C}^{\infty}_0\biggl([0,\eps),\begin{array}{c} \dot{C}^{\infty}(H,F) \\ \oplus \\ \dot{C}^{\infty}(\regbdry H,J_+) \end{array}\biggr).
$$
In coordinates $t = -\frac{1}{x}$, this operator is of the form \eqref{transformedoperator} with a global symbol
$b(t,t',\tau,\lambda)$ in the class
\begin{equation}\label{S00psicu}
S_{\cl}^{0,0}(\R_t\times\R_{t'},\Psicu^{\mu,d}(H,\R_{\tau}\times\Lambda)) =
S_{\cl}^0(\R_t) \hat{\otimes}_{\pi} S_{\cl}^0(\R_{t'}) \hat{\otimes}_{\pi} \Psicu^{\mu,d}(H,\R_{\tau}\times\Lambda)
\end{equation}
supported in $\{(t,t') \st t,t' < 0\}$.
By induction, $\Psicu^{\mu,d}(H,\R_{\tau}\times\Lambda)$ embeds into a space of $(\tau,\lambda)$-dependent
operator-valued symbols between the cusp Sobolev spaces on $H$. Consequently, operator-valued variants of
oscillatory integral techniques, which were introduced in \cite{Cordes,Parenti} for scalar symbols on
Euclidean space, are applicable here. Such arguments, also known as \emph{Kumano-go's technique}, have been employed at
various occasions in the literature on pseudodifferential operators. We use them implicitly throughout
this section for setting up the asymptotic properties of our calculus (see also \cite{KraiParabBVP}, which relates
to manifolds with cylindrical ends).
The explicit iterative construction of $\Psicu^{\mu,d}(H,\R_{\tau}\times\Lambda)$ shows that
the symbol classes \eqref{S00psicu} remain indeed preserved under all manipulations in Kumano-go's technique that are
involved in this process.
\end{remark}

\begin{definition}\label{properlysupported}
Formally, we write the operator $A(\lambda) \in \Psicu^{\mu,d}([0,\eps)\times H,\Lambda)$ from
\eqref{cuspquantization} as
$$
\bigl(A(\lambda)u\bigr)(x) = \int\limits_{0}^{\eps}K_A(x,y,\lambda)u(y)\,\frac{dy}{y^2}
$$
with the operator-valued distributional kernel
$$
K_A(x,y,\lambda) = \frac{1}{2\pi}\int\limits_{\R}e^{i(1/y-1/x)\xi}a(x,\xi,\lambda)\,d\xi + c(x,y,\lambda)
$$
in the variables $(x,y) \in (0,\eps)^2$ and parameter $\lambda \in \Lambda$.

We call $A(\lambda)$ \emph{properly supported up to the origin}, if for each $0 < \delta < \eps$ there exist compact
sets $K_x \subset [0,\delta]\times[0,\eps)$ and $K_y \subset [0,\eps)\times[0,\delta]$ independent
of $\lambda \in \Lambda$ such that
\begin{align*}
\overline{\supp K_A(\lambda)} \cap \bigl([0,\delta]\times[0,\eps)\bigr) &\subset K_x, \\
\intertext{and}
\overline{\supp K_A(\lambda)} \cap \bigl([0,\eps)\times[0,\delta]\bigr) &\subset K_y,
\end{align*}
where $\overline{\supp K_A(\lambda)}$ is the closure of $\supp K_A(\lambda)$ in $[0,\eps)^2$.

Observe that this notion of properly supportedness up to the origin is \emph{not compatible} with the usual notion
of properly supportedness of pseudodifferential operators on a smooth boundaryless manifold. However, if
$K_A(x,y,\lambda) \equiv 0$ for $\{x < \delta\}\ \cup \{y < \delta\}$, then $A(\lambda)$ is properly supported
up to the origin if and only if $A(\lambda)$ is properly supported in the usual (operator-valued) sense on $(0,\eps)$.
\end{definition}

\begin{proposition}\label{Propsuppmapprop}
Let $A(\lambda) \in \Psicu^{\mu,d}([0,\eps)\times H,\Lambda)$ be properly supported up to the origin. Then $A(\lambda)$
is continuous in the spaces
\begin{align*}
A(\lambda) &: \dot{C}^{\infty}_0\biggl([0,\eps),\begin{array}{c} \dot{C}^{\infty}(H,E) \\ \oplus \\ \dot{C}^{\infty}(\regbdry H,J_-) \end{array}\biggr) \to
\dot{C}_0^{\infty}\biggl([0,\eps),\begin{array}{c} \dot{C}^{\infty}(H,F) \\ \oplus \\ \dot{C}^{\infty}(\regbdry H,J_+) \end{array}\biggr), \\
A(\lambda) &: \dot{C}^{\infty}\biggl([0,\eps),\begin{array}{c} \dot{C}^{\infty}(H,E) \\ \oplus \\ \dot{C}^{\infty}(\regbdry H,J_-) \end{array}\biggr) \to
\dot{C}^{\infty}\biggl([0,\eps),\begin{array}{c} \dot{C}^{\infty}(H,F) \\ \oplus \\ \dot{C}^{\infty}(\regbdry H,J_+) \end{array}\biggr).
\end{align*}
\end{proposition}
\begin{proof}
Pick cut-off functions $\hat{\omega} \prec \omega \prec \tilde{\omega}$ near zero, and write
$$
A(\lambda) = \omega A(\lambda)\tilde{\omega} + (1-\omega)A(\lambda)(1-\hat{\omega}) + R(\lambda).
$$
Each summand is properly supported up to the origin. By Lemma \ref{localpseudolocal} the term $R(\lambda)$ is of the
form \eqref{Cdef}, and thus trivially has the asserted mapping properties. Obviously, this is also the case for
$\omega A(\lambda) \tilde{\omega}$. $(1-\omega)A(\lambda)(1-\hat{\omega})$ is properly supported in the usual
sense on $(0,\eps)$, and its operator-valued Schwartz kernel is supported away from $\{x=0\}\cup\{y=0\}$.
Consequenty, also this term has the desired mapping properties.
\end{proof}

\begin{lemma}\label{properlysupportedrepr}
Every $A(\lambda) \in \Psicu^{\mu,d}([0,\eps)\times H,\Lambda)$ can be written
in the form \eqref{cuspquantization}
$$
\bigl(A(\lambda)u\bigr)(x) = \frac{1}{2\pi}\int\limits_{\R}\int\limits_{0}^{\eps}
e^{i(1/y-1/x)\xi}a(x,\xi,\lambda)u(y)\,\frac{dy}{y^2}\,d\xi + \bigl(C(\lambda)u\bigr)(x)
$$
with $A(\lambda) - C(\lambda)$ properly supported up to the origin.
\end{lemma}
\begin{proof}
Let $\hat{\omega} \prec \omega \prec \tilde{\omega}$ be cut-off functions near zero, and write
$$
A(\lambda) = \omega A(\lambda)\tilde{\omega} + (1-\omega)A(\lambda)(1-\hat{\omega}) + R(\lambda).
$$
$R(\lambda)$ is of the form \eqref{Cdef} by Lemma \ref{localpseudolocal}, and the operator
$\omega A(\lambda)\tilde{\omega}$ is properly supported up to the origin and has a representation of the form
\begin{equation}\label{omegaAomega}
\bigl(\omega A(\lambda) \tilde{\omega} u\bigr)(x) = \frac{1}{2\pi}\int\limits_{\R}\int\limits_{0}^{\eps}
e^{i(1/y-1/x)\xi}\tilde{a}(x,\xi,\lambda)u(y)\,\frac{dy}{y^2}\,d\xi
\end{equation}
with a unique symbol $\tilde{a}(x,\xi,\lambda)$, because the operator-valued Schwartz kernel of
$\omega A(\lambda) \tilde{\omega}$ is supported in $[0,\delta]^2$ for some $0 < \delta < \eps$
(see Remark \ref{transformation}).

For the analyis of the term $(1-\omega)A(\lambda)(1-\hat{\omega})$ write
$A(\lambda) = A_{\textup{prop}}(\lambda) + C'(\lambda)$, where
\begin{align}
\bigl(C'(\lambda)u\bigr)u\bigr)(x) &= \int\limits_0^{\eps}c'(x,y,\lambda)u(y)\,\frac{dy}{y^2}, \\
c'(x,y,\lambda) &\in C^{\infty}((0,\eps)_x\times(0,\eps)_y,\Psicu^{-\infty,d}(H,\Lambda)), \\
\intertext{and}
\label{Aprop}\bigl(A_{\textup{prop}}(\lambda)u\bigr)(x) &= \frac{1}{2\pi}\iint e^{i(1/y-1/x)\xi}a_{\textup{prop}}(x,\xi,\lambda)u(y)\,\frac{dy}{y^2}\,d\xi, \\
a_{\textup{prop}}(x,\xi,\lambda) &\in C^{\infty}((0,\eps)_x,\Psicu^{\mu,d}(H,\R_{\xi}\times\Lambda)),
\end{align}
and $A_{\textup{prop}}(\lambda)$ is properly supported on $(0,\eps)$ in the usual sense.
Note that $(1-\omega)C'(\lambda)(1-\hat{\omega})$ is of the form \eqref{Cdef}, while
$(1-\omega)A_{\textup{prop}}(\lambda)(1-\hat{\omega})$ is properly supported on $(0,\eps)$, and thus can be written
in the form \eqref{Aprop} with some (other) symbol $a'_{\textup{prop}}(x,\xi,\lambda)$ which vanishes for small $x$.

Hence the assertion follows with $a = \tilde{a} + a'_{\textup{prop}}$.
\end{proof}

\begin{definition}\label{localpseudoconormal}
Let $A(\lambda) \in \Psicu^{\mu,d}([0,\eps)\times H,\Lambda)$ be written according to \eqref{cuspquantization}
with a symbol $a(x,\xi,\lambda)$ as in \eqref{cuspsymbol}.
The \emph{conormal symbol} (or normal operator) of $A(\lambda)$ with respect to $H$ is defined as
\begin{equation}\label{localconormaldef}
N_H(A)(\xi,\lambda):= a(0,\xi,\lambda) \in \Psicu^{\mu,d}(H,\R\times\Lambda),
\end{equation}
and it can be regarded as a family of operators
\begin{equation}\label{conormalsymbolcinf}
N_H(A)(\xi,\lambda) : \begin{array}{c} \dot{C}^{\infty}(H,E) \\ \oplus \\ \dot{C}^{\infty}(\regbdry H,J_-) \end{array}
\to \begin{array}{c} \dot{C}^{\infty}(H,F) \\ \oplus \\ \dot{C}^{\infty}(\regbdry H,J_+) \end{array}
\end{equation}
or, alternatively, as a family of operators in the cusp Sobolev spaces on $H$ and $\regbdry H$.
\end{definition}

The conormal symbol is indeed well defined for the operator $A(\lambda)$,
i.e. independent of the choice of the symbol $a(x,\xi,\lambda)$ that is involved in the representation
\eqref{cuspquantization} of $A(\lambda)$. This is due to the following:

Let $\hat{\omega} \prec \omega \prec \tilde{\omega}$ be cut-off functions near zero, and write
$$
A(\lambda) = \omega A(\lambda)\tilde{\omega} + (1-\omega)A(\lambda)(1-\hat{\omega}) + R(\lambda).
$$
$R(\lambda)$ is of the form \eqref{Cdef}, and the term $(1-\omega)A(\lambda)(1-\hat{\omega})$
clearly also does not contribute to $N_H(A)$ because its operator-valued Schwartz kernel is supported strictly away
from $\{x=0\}\cup\{y=0\}$. The operator $\omega A(\lambda) \tilde{\omega}$ has a representation of the
form \eqref{omegaAomega} with a unique symbol $\tilde{a}(x,\xi,\lambda)$, and
$a(x,\xi,\lambda) - \tilde{a}(x,\xi,\lambda)$ vanishes to infinite order at $x = 0$.
In particular, $a(0,\xi,\lambda) = \tilde{a}(0,\xi,\lambda)$ does not depend on the specific representative $a$.

\begin{proposition}\label{localcomposition}
Let $A_j(\lambda) \in \Psicu^{\mu_j,d_j}([0,\eps)\times H,\Lambda)$, $j = 1,2$, and either $A_1(\lambda)$ or
$A_2(\lambda)$ properly supported up to the origin such that the composition $A_1(\lambda)A_2(\lambda)$ is well defined
according to Proposition \ref{Propsuppmapprop} (the vector bundles are assumed to fit together).

Then $A_1(\lambda)A_2(\lambda) \in \Psicu^{\mu_1+\mu_2,\max\{d_1+\mu_2,d_2\}}([0,\eps)\times H,\Lambda)$. More
precisely, if
$$
\bigl(A_j(\lambda)u\bigr)(x) = \frac{1}{2\pi}\int\limits_{\R}\int\limits_{0}^{\eps}
e^{i(1/y-1/x)\xi}a_j(x,\xi,\lambda)u(y)\,\frac{dy}{y^2}\,d\xi + \bigl(C_j(\lambda)u\bigr)(x),
$$
$j=1,2$, are representations according to \eqref{cuspquantization}, then
$$
\bigl(A_1(\lambda)A_2(\lambda)u\bigr)(x) = \frac{1}{2\pi}\int\limits_{\R}\int\limits_{0}^{\eps}
e^{i(1/y-1/x)\xi}(a_1\# a_2)(x,\xi,\lambda)u(y)\,\frac{dy}{y^2}\,d\xi + \bigl(C(\lambda)u\bigr)(x)
$$
is a representation of the composition with a symbol
$$
(a_1\# a_2)(x,\xi,\lambda) \in C^{\infty}\bigl([0,\eps)_x,\Psicu^{\mu_1+\mu_2,\max\{d_1+\mu_2,d_2\}}(H,\R_{\xi}\times\Lambda)\bigr)
$$
that has an asymptotic expansion
\begin{equation}\label{Leibnizasymp}
a_1\# a_2 \sim \sum\limits_{k=0}^{\infty}\frac{1}{k!}\bigl(\partial^k_{\xi}a_1\bigr)\bigl((x^2D_x)^ka_2\bigr)
\end{equation}
in the sense that the difference
$a_1\# a_2 - \sum\limits_{k=0}^{K-1}\frac{1}{k!}\bigl(\partial^k_{\xi}a_1\bigr)\bigl((x^2D_x)^ka_2\bigr)$
belongs to the space
$x^K C^{\infty}\bigl([0,\eps)_x,\Psicu^{\mu_1+\mu_2-K,\max\{d_1+\mu_2,d_2\}}(H,\R_{\xi}\times\Lambda)\bigr)$
for every $K \in \N_0$.

In particular, we have
\begin{equation}\label{princsymbolcomp}
\left.\begin{aligned}
\cusym(A_1A_2) &= \cusym(A_1)\cusym(A_2), \\
\cusymb(A_1A_2) &= \cusymb(A_1)\cusymb(A_2), \\
N_H(A_1A_2) &= N_H(A_1)N_H(A_2)
\end{aligned}\right\}
\end{equation}
for the cusp-principal symbols and the conormal symbol of the composition.
\end{proposition}
\begin{proof}
Let $\omega,\tilde{\omega} \in C_0^{\infty}([0,\eps))$ be cut-off functions near zero. As either $A_1(\lambda)$ or
$A_2(\lambda)$ is properly supported up to the origin, there exist cut-off functions $\hat{\omega},\check{\omega} \in C_0^{\infty}([0,\eps))$
near zero with $\hat{\omega} \prec \check{\omega}$ such that
\begin{equation}\label{omegaA1A2omega}
\omega A_1(\lambda) A_2(\lambda) \tilde{\omega} \equiv \bigl(\omega A_1(\lambda) \hat{\omega}\bigr)\bigl(\check{\omega} A_2(\lambda) \tilde{\omega}\bigr).
\end{equation}
Hence, if any of the $A_j(\lambda)$ is of the form \eqref{Cdef}, we conclude that also $\omega A_1(\lambda) A_2(\lambda) \tilde{\omega}$
is of the form \eqref{Cdef} by transforming \eqref{omegaA1A2omega} to the coordinate $t = -\frac{1}{x}$ and employing
Kumano-go's technique, see Remark \ref{transformation}. Consequently also $A_1(\lambda)A_2(\lambda)$ is
of the form \eqref{Cdef}, i.e. the integral operators \eqref{Cdef} with kernels \eqref{opkernel} form a two-sided
ideal.

Now let $\omega \prec \hat{\omega} \prec \check{\omega} \prec \tilde{\omega}$. By Lemma \ref{localpseudolocal} and what
we just proved, we conclude that the equality \eqref{omegaA1A2omega} holds for this choice of cut-off functions
modulo an integral operator of the form \eqref{Cdef} (with appropriate type).
Using again a change of variables, Kumano-go's technique, and the induction hypothesis as regards the parameter-dependent
cusp calculus of boundary value problems on $H$, we conclude that the right-hand side of \eqref{omegaA1A2omega} is
of the form \eqref{cuspquantization} with $C(\lambda) \equiv 0$ and an operator-valued symbol \eqref{cuspsymbol} of order
$\mu_1+\mu_2$ and type $\max\{d_1+\mu_2,d_2\}$, which behaves asymptotically like $\omega a_1 \# a_2$ in the sense specified by
\eqref{Leibnizasymp}.

Consider now cut-off functions $\tilde{\omega} \prec \check{\omega} \prec \hat{\omega} \prec \omega$ near zero, and the
operator $(1-\omega)A_1(\lambda)A_2(\lambda)(1-\tilde{\omega})$. By Lemma \ref{localpseudolocal} and the ideal
property of the operators \eqref{Cdef} we conclude that
\begin{equation}\label{1moA1mo}
(1-\omega)A_1(\lambda)A_2(\lambda)(1-\tilde{\omega}) \equiv \bigl((1-\omega)A_1(\lambda)(1-\hat{\omega})\bigr)
\bigl((1-\check{\omega})A_2(\lambda)(1-\tilde{\omega})\bigr)
\end{equation}
modulo an integral operator of the form \eqref{Cdef}. Writing $A_j(\lambda) = A_{j,\textup{prop}}(\lambda) + C_j'(\lambda)$,
$j=1,2$, as in the proof of Lemma \ref{properlysupportedrepr}, and using again the ideal property of the operators \eqref{Cdef},
we see that we only have to analyze the composition on the right-hand side of \eqref{1moA1mo} where, in addition, the $A_j(\lambda)$ can both be replaced by
their properly supported representatives $A_{j,\textup{prop}}(\lambda)$. This composition, however, can be handled
with ordinary techniques from the theory of pseudodifferential operators with operator-valued symbols, and from the theory
of pseudodifferential boundary value problems --- that this is indeed the case follows from our induction hypothesis resp. the iterative construction of the cusp
calculus of pseudodifferential boundary value problems. As a result, we obtain that the composition \eqref{1moA1mo}
is of the form \eqref{cuspquantization} with an operator-valued symbol \eqref{cuspsymbol} of order $\mu_1+\mu_2$ and type
$\max\{d_1+\mu_2,d_2\}$ which behaves asymptotically like $(1-\omega) a_1 \# a_2$, see \eqref{Leibnizasymp}.

Finally, if $\omega \prec \tilde{\omega}$, then it is immediate from Lemma \ref{localpseudolocal} and the first part
of this proof that both compositions
$\omega A_1(\lambda)A_2(\lambda) (1-\tilde{\omega})$ and $(1-\tilde{\omega})A_1(\lambda)A_2(\lambda)\omega$ are
of the form \eqref{Cdef}. This completes the proof of the proposition.
\end{proof}

\begin{remark}\label{vanishinglocal}
Let $A(\lambda) \in \Psicu^{\mu,d}([0,\eps)\times H,\Lambda)$ and $K \in \Z$. Following the proof of
Proposition \ref{localcomposition} we obtain that
$$
B(\lambda) = x^KA(\lambda)x^{-K} \in \Psicu^{\mu,d}([0,\eps)\times H,\Lambda),
$$
and $\cusym(A) = \cusym(B)$, $\cusymb(A) = \cusymb(B)$, as well as $N_H(A) = N_H(B)$.

This shows, in particular, that the spaces $x^K\,\Psicu^{\ast,\ast}([0,\eps)\times H,\Lambda)$ form two-sided ideals
in $\Psicu^{\ast,\ast}([0,\eps)\times H,\Lambda)$ for every $K \in \N_0$.
\end{remark}

\begin{definition}\label{localcuspelliptic}
A boundary value problem $A(\lambda) \in \Psicu^{\mu,d}([0,\eps)\times H,\Lambda)$ is called \emph{cusp-elliptic} (with parameter
$\lambda \in \Lambda$), if both the cusp-principal symbol $\cusym(A)$ and the cusp-principal boundary symbol
$\cusymb(A)$ are invertible on $\bigl(\cuT^*\bigl([0,\eps)\times H\bigr) \times \Lambda\bigr)\setminus 0$
and $\bigl(\cuT^*\bigl([0,\eps)\times \regbdry H\bigr) \times \Lambda\bigr) \setminus 0$, respectively.
\end{definition}

\begin{proposition}\label{localcuspparametrix}
Let $A(\lambda) \in \Psicu^{\mu,d}([0,\eps)\times H,\Lambda)$ be cusp-elliptic. Then there exists a parametrix
$B(\lambda) \in \Psicu^{-\mu,(d-\mu)_+}([0,\eps)\times H,\Lambda)$, where $(d-\mu)_+ = \max\{0,d-\mu\}$, such that
\begin{align*}
B(\lambda)A(\lambda) - 1 &\in \Psicu^{-\infty,\ast}([0,\eps)\times H,\Lambda), \quad \text{and} \\
A(\lambda)B(\lambda) - 1 &\in \Psicu^{-\infty,\ast}([0,\eps)\times H,\Lambda).
\end{align*}
The types of these remainders are given by the type formula for the composition from Proposition \ref{localcomposition}.

If, in addition, the conormal symbol $N_H(A)(\xi,\lambda)$ is invertible in the spaces \eqref{conormalsymbolcinf} for
all $(\xi,\lambda) \in \R\times\Lambda$, then there exists $R(\lambda) \in \Psicu^{-\infty,(d-\mu)_+}([0,\eps)\times H,\Lambda)$
such that both $\bigl(B(\lambda)+R(\lambda)\bigr)A(\lambda) - 1$ and $A(\lambda)\bigl(B(\lambda)+R(\lambda)\bigr)-1$
are of the form \eqref{Cdef} with kernels \eqref{opkernel} (with the appropriate types).

Both $B(\lambda)$ and $R(\lambda)$ are properly supported up to the origin.
\end{proposition}
\begin{proof}
Let $a(x,\xi,\lambda)$ be the symbol in a representation \eqref{cuspquantization} for $A(\lambda)$.
The cusp-ellipticity of $A(\lambda)$ implies the existence of
\begin{align*}
b'(x,\xi,\lambda) &\in C^{\infty}\bigl([0,\eps)_x,\Psicu^{-\mu,(d-\mu)_+}(H,\R_{\xi}\times\Lambda)\bigr) \\
\intertext{such that}
\left.\begin{aligned}
a(x,\xi,\lambda)b'(x,\xi,\lambda)&-1 \\
b'(x,\xi,\lambda)a(x,\xi,\lambda)&-1
\end{aligned}\right\} &\in C^{\infty}\bigl([0,\eps)_x,\Psicu^{-\infty,\ast}(H,\R_{\xi}\times\Lambda)\bigr),
\end{align*}
see Theorem \ref{Parametrix} -- this argument makes use of our induction hypothesis.
By Lemma \ref{properlysupportedrepr}, Proposition \ref{localcomposition}, and a
standard formal Neumann series argument -- which is applicable here by the induction hypothesis -- we
conclude that there exists a symbol $b(x,\xi,\lambda)$ such that
$$
\left.\begin{aligned}
a(x,\xi,\lambda)\# b(x,\xi,\lambda)&-1 \\
b(x,\xi,\lambda)\# a(x,\xi,\lambda)&-1
\end{aligned}\right\} \in C^{\infty}\bigl([0,\eps)_x,\Psicu^{-\infty,\ast}(H,\R_{\xi}\times\Lambda)\bigr).
$$
Hence the assertion of the proposition follows with $B(\lambda)$ given by
$$
\bigl(B(\lambda)u\bigr)(x) = \frac{1}{2\pi}\int\limits_{\R}\int\limits_{0}^{\eps}
e^{i(1/y-1/x)\xi}b(x,\xi,\lambda)u(y)\,\frac{dy}{y^2}\,d\xi,
$$
and $B(\lambda)$ is properly supported up to the origin.

Let us now assume that, in addition, $N_H(A)(\xi,\lambda) = a(0,\xi,\lambda)$ is invertible for all $(\xi,\lambda) \in \R\times\Lambda$.
By induction, we conclude that the inverse $a(0,\xi,\lambda)^{-1}$ belongs to the space $\Psicu^{-\mu,(d-\mu)_+}(H,\R\times\Lambda)$, see
Theorem \ref{largeparaminverse}. By the multiplicativity of the conormal symbols, see Proposition \ref{localcomposition},
we further conclude that
$$
N_H(A)(\xi,\lambda)^{-1} - N_H(B)(\xi,\lambda) = r'(\xi,\lambda) \in \Psicu^{-\infty,(d-\mu)_+}(H,\R\times\Lambda),
$$
and by quantizing this operator-valued symbol according to \eqref{cuspquantization}, we see that
$r'(\xi,\lambda) = N_H(R')(\xi,\lambda)$ for some $R'(\lambda) \in \Psicu^{-\infty,(d-\mu)_+}([0,\eps)\times H,\Lambda)$.
Here we may assume that $R'(\lambda)$ is properly supported up to the origin, otherwise we substitute $R'(\lambda)$ by
$\omega R'(\lambda) \omega$, where $\omega \in C_0^{\infty}([0,\eps))$ is any cut-off function near zero. Consequently,
$$
\left.\begin{aligned}
\bigl(B(\lambda)+R'(\lambda)\bigr)A(\lambda) &- 1 \\
A(\lambda)\bigl(B(\lambda)+R'(\lambda)\bigr) &- 1
\end{aligned}
\right\} \in x\Psicu^{-\infty,\ast}([0,\eps)\times H,\Lambda),
$$
and thus it remains to show that an operator of the form $1 + G(\lambda)$
with $G(\lambda) \in x\,\Psicu^{-\infty,d}([0,\eps)\times H,\Lambda)$ has a parametix $1 + G'(\lambda)$
up to remainders of the form \eqref{Cdef}, where $G'(\lambda) \in x\,\Psicu^{-\infty,d}([0,\eps)\times H,\Lambda)$.

For this proof we may assume by Lemma \ref{properlysupportedrepr} that $G(\lambda)$ is properly supported up to the origin.
Let
$$
g_k(x,\xi,\lambda) \in x^kC^{\infty}\bigl([0,\eps)_x,\Psicu^{-\infty,d}(H,\R_{\xi}\times\Lambda)\bigr)
$$
be an operator-valued symbol \eqref{cuspsymbol} in the representation \eqref{cuspquantization} for the composition $G(\lambda)^k \in x^k\,\Psicu^{-\infty,d}([0,\eps)\times H,\Lambda)$, $k \in \N$.
A Borel argument shows the existence of a symbol $g'(x,\xi,\lambda)$ such that
$$
g'(x,\xi,\lambda) - \sum\limits_{k=1}^{K-1}(-1)^kg_k(x,\xi,\lambda) \in x^KC^{\infty}\bigl([0,\eps)_x,\Psicu^{-\infty,d}(H,\R_{\xi}\times\Lambda)\bigr)
$$
for $K \in \N$. Now pick $D(\lambda)$ of the form \eqref{Cdef} such that
$$
\bigl(G'(\lambda)u\bigr)(x):= \frac{1}{2\pi}\int\limits_{\R}\int\limits_{0}^{\eps}
e^{i(1/y-1/x)\xi}g'(x,\xi,\lambda)u(y)\,\frac{dy}{y^2}\,d\xi + \bigl(D(\lambda)u\bigr)(x)
$$
is properly supported up to the origin. Then $G'(\lambda) \in x\,\Psicu^{-\infty,d}([0,\eps)\times H,\Lambda)$, and
by construction we have
$$
\left.\begin{aligned}
\bigl(1+G(\lambda)\bigr)\bigl(1+G'(\lambda)\bigr) &- 1 \\
\bigl(1+G'(\lambda)\bigr)\bigl(1+G(\lambda)\bigr) &- 1
\end{aligned}\right\}
\in \bigcap\limits_{K \in \N}x^K\,\Psicu^{-\infty,d}([0,\eps)\times H,\Lambda),
$$
where the latter is just the space of all integral operators \eqref{Cdef} with kernels \eqref{opkernel}.
This completes the proof of the proposition.
\end{proof}

\subsection{The full cusp algebra on $\Mbar$}\label{subsec-Fullalgebra}

This section is devoted to set up the class $\Psicu^{\mu,d}(\Mbar,\Lambda)$ of (parameter-dependent)
pseudodifferential boundary value problems in the cusp algebra. The operators
\begin{equation}\label{AlC0}
A(\lambda) : \begin{array}{c} C_0^{\infty}(\overline{N},E) \\ \oplus \\ C_0^{\infty}(\partial\overline{N},J_-) \end{array} \to
\begin{array}{c} C^{\infty}(\overline{N},F) \\ \oplus \\ C^{\infty}(\partial\overline{N},J_+) \end{array}
\end{equation}
in $\Psicu^{\mu,d}(\Mbar,\Lambda)$ belong to the class $\B^{\mu,d}(\overline{N},\Lambda)$
of (parameter-dependent) operators in Boutet de Monvel's calculus on $\overline{N} = \Mbar\setminus\singbdry\Mbar$
with a specific behavior near $\singbdry\Mbar$:

\begin{definition}\label{cuspalgebradef}
Let $\mu \in \Z$ and $d \in \N_0$. An operator $A(\lambda)$ in the spaces \eqref{AlC0} belongs to the class
$\Psicu^{\mu,d}(\Mbar,\Lambda)$ if and only if the following holds (we employ here throughout the standard
convention \eqref{Multop} for the multiplication operator with a function $\varphi$):
\begin{enumerate}[i)]
\item Let $\varphi_1,\varphi_2 \in C^{\infty}(\Mbar)$ with $\supp \varphi_1 \cap \supp \varphi_2 = \emptyset$, and assume that
each singular hypersurface $H \subset \singbdry\Mbar$ has nonempty intersection with at most one of the supports of the $\varphi_j$'s.
Then the operator $\varphi_1 A(\lambda) \varphi_2$ is required to belong to the space $\Psicu^{-\infty,d}_{\infty}(\Mbar,\Lambda)$ of
residual Green operators of type $d$, see Section \ref{subsec-ResidualGreen}.
\item Let $\varphi_1,\varphi_2 \in C^{\infty}(\Mbar)$ be such that $\supp\varphi_j \cap \singbdry\Mbar = \emptyset$
for $j = 1,2$. Then we require $\varphi_1 A(\lambda) \varphi_2$ to belong to the class $\B^{\mu,d}(\overline{N},\Lambda)$ of
pseudodifferential boundary value problems in Boutet de Monvel's calculus on $\overline{N}$.
\item Let $\varphi_1,\varphi_2 \in C^{\infty}(\Mbar)$ be supported inside the same collar neighborhood
$[0,\eps)\times H$ of a singular boundary hypersurface $H \subset \singbdry\Mbar$, and assume that
$H \cap \regbdry\Mbar \neq \emptyset$. Then $\varphi_1 A(\lambda) \varphi_2$ is required to belong to the class
$\Psicu^{\mu,d}([0,\eps)\times H,\Lambda)$ discussed in Section \ref{subsec-SingularCollar}.
\item Let $\varphi_1,\varphi_2 \in C^{\infty}(\Mbar)$ be supported inside the same collar neighborhood
$[0,\eps)\times H$ of a singular boundary hypersurface $H \subset \singbdry\Mbar$, and assume now that
$H \cap \regbdry\Mbar = \emptyset$. Then the only nonzero term in the matrix operator $\varphi_1 A(\lambda) \varphi_2$
is the interior operator in the upper left corner, which is required to belong to the class of ordinary
(parameter-dependent) cusp operators $\Psicu^{\mu}([0,\eps)\times H,\Lambda)$ near the hypersurface $H$ (see the
notes at the beginning of Section \ref{subsec-SingularCollar}).
\end{enumerate}
\end{definition}
The projective topology with respect to the mappings
$$
\Psicu^{\mu,d}(\Mbar,\Lambda) \ni A(\lambda) \mapsto \varphi_1 A(\lambda) \varphi_2 \in
\begin{cases}
\Psicu^{-\infty,d}_{\infty}(\Mbar,\Lambda) \\
\B^{\mu,d}(\overline{N},\Lambda) \\
\Psicu^{\mu(,d)}([0,\eps)\times H,\Lambda)
\end{cases}
$$
according to i)--iv) in Definition \ref{cuspalgebradef} makes $\Psicu^{\mu,d}(\Mbar,\Lambda)$ a Fr{\'e}chet space.

It is evident from Definition \ref{cuspalgebradef} that every $A(\lambda) \in \Psicu^{\mu,d}(\Mbar,\Lambda)$ extends to
$$
A(\lambda) : \begin{array}{c} \dot{C}^{\infty}(\Mbar,E) \\ \oplus \\ \dot{C}^{\infty}(\regbdry\Mbar,J_-) \end{array}
\to \begin{array}{c} \dot{C}^{\infty}(\Mbar,F) \\ \oplus \\ \dot{C}^{\infty}(\regbdry\Mbar,J_+) \end{array},
$$
see \eqref{initialcusp}. We will henceforth consider these spaces a core for the operators in the cusp algebra
of boundary value problems. Moreover, by the iterative construction of the calculus, we obtain the following

\begin{proposition}\label{extensionbycontinuity}
Every $A(\lambda) \in \Psicu^{\mu,d}(\Mbar,\Lambda)$ extends by continuity to
\begin{equation}\label{extensionSobCont}
A(\lambda) : \begin{array}{c} \x_{\sing}^{\alpha}\Hcu^s(\Mbar,E) \\ \oplus \\ \x_{\sing}^{\alpha}\Hcu^s(\regbdry\Mbar,J_-) \end{array} \to
\begin{array}{c} \x_{\sing}^{\alpha}\Hcu^{s-\mu}(\Mbar,F) \\ \oplus \\ \x_{\sing}^{\alpha}\Hcu^{s-\mu}(\regbdry\Mbar,J_+) \end{array}
\end{equation}
for $s > d-\frac{1}{2}$ and all $\alpha \in \R^{N-\ell}$, and the class $\Psicu^{\ast,d}(\Mbar,\Lambda)$ embeds into
the space of symbols depending on the parameter $\lambda \in \Lambda$ that are operator-valued in the bounded operators
between the cusp Sobolev spaces.
Recall that $\singbdry\Mbar$ consists of $N-\ell$, and $\regbdry\Mbar$ of $\ell$ hypersurfaces.
\end{proposition}

\begin{definition}\label{principalsymbolcusppseudo}
The operators $A(\lambda) = \bigl(A_{i,j}(\lambda)\bigr)_{i,j=1,2} \in \Psicu^{\mu,d}(\Mbar,\Lambda)$ have the following principal symbols (see
Definition \ref{principalsymbolscuspdifferential}):
\begin{enumerate}[i)]
\item The \emph{cusp-principal symbol}
$$
\cusym(A) \in C^{\infty}\bigl(\bigl(\cuT^*\Mbar\times\Lambda\bigr)\setminus 0,\Hom(\cupi^*E,\cupi^*F)\bigr),
$$
which is the canonical extension to the cusp cotangent bundle $\bigl(\cuT^*\Mbar\times\Lambda\bigr)\setminus 0$
of the homogeneous principal symbol $\sym(A_{1,1})$ of the (parameter-dependent) pseudodifferential
operator $A_{1,1}(\lambda)$ in the upper left corner of the operator matrix $A(\lambda)$ (see also Section \ref{subsec-SingularCollar}).
Note that the principal symbol $\sym(A_{1,1})$ is defined initially only on $\bigl(T^*\overline{N}\times\Lambda\bigr)\setminus 0$.

The cusp-principal symbol $\cusym(A)$ is a homogeneous function of degree $\mu$ in the fibres of
$\bigl(\cuT^*\Mbar\times\Lambda\bigr)\setminus 0$.
\item The \emph{cusp-principal boundary symbol} $\cusymb(A)$, a section in
$$
C^{\infty}\biggl(\bigl(\cuT^*\regbdry\Mbar\times\Lambda\bigr)\setminus 0,
\Hom\biggl(
\begin{array}{c} \Scu\otimes\cupi^*E|_{\regbdry\Mbar} \\ \oplus \\ \cupi^*J_- \end{array},
\begin{array}{c} \Scu\otimes\cupi^*F|_{\regbdry\Mbar} \\ \oplus \\ \cupi^*J_+ \end{array}
\biggr)\biggr).
$$
Similar to the cusp-principal symbol, the cusp-principal boundary symbol is the canonical extension of the
principal boundary symbol
$$
\sym_{\partial}(A) \in
C^{\infty}\biggl(\bigl(T^*\partial\overline{N}\times\Lambda\bigr)\setminus 0,
\Hom\biggl(
\begin{array}{c} \S\otimes\pi^*E|_{\partial\overline{N}} \\ \oplus \\ \pi^*J_- \end{array},
\begin{array}{c} \S\otimes\pi^*F|_{\partial\overline{N}} \\ \oplus \\ \pi^*J_+ \end{array}
\biggr)\biggr)
$$
of the operator $A(\lambda) \in \B^{\mu,d}(\overline{N},\Lambda)$ to $\bigl(\cuT^*\regbdry\Mbar\times\Lambda\bigr)\setminus 0$,
see also the discussion around \eqref{cuspbdrysymbol} and in Section \ref{subsec-SingularCollar}.

The cusp-principal boundary symbol $\cusymb(A)$ is $\kappa$-homogeneous of degree $\mu$, i.e.
$\cusymb(A)(z,\varrho\zeta,\varrho\lambda)$ equals
$$
\varrho^{\mu}
\begin{pmatrix}\kappa_{\varrho}\otimes\id_{\cupi^*F|_{\regbdry\Mbar}} & 0 \\ 0 & \id_{\cupi^*J_+} \end{pmatrix}
\cusymb(A)(z,\zeta,\lambda)
\begin{pmatrix}\kappa^{-1}_{\varrho}\otimes\id_{\cupi^*E|_{\regbdry\Mbar}} & 0 \\ 0 & \id_{\cupi^*J_-} \end{pmatrix}
$$
for $\varrho > 0$ and $(z,\zeta,\lambda) \in \bigl(\cuT^*\regbdry\Mbar\times\Lambda\bigr)\setminus 0$. Here
$\kappa_{\varrho}$ is the natural $\R_+$-action \eqref{kappagroup} in the $\S(\overline{\R}_+)$-fibres of the (lifted) bundle $\Scu \to
\bigl(\cuT^*\regbdry\Mbar\times\Lambda\bigr)\setminus 0$.
\end{enumerate}
Let ${}^{\cu}\Sigma \equiv {}^{\cu}\Sigma(\Mbar,\Lambda)$ be the space of cusp-principal symbols of the operators in
$\Psicu^{\mu,d}(\Mbar,\Lambda)$, i.e. the space of tuples of homogeneous and $\kappa$-homogeneous
sections that satisfy a canonical compatibility condition. By means of a partition of unity and the local splittings of
the principal symbol sequences in $\B^{\mu,d}(\overline{N},\Lambda)$ and $\Psicu^{\mu(,d)}([0,\eps)\times H,\Lambda)$,
see Proposition \ref{cusplocalsplitexact}, we obtain that the cusp-principal symbol sequence
$$
\begin{CD}
0 @>>> \Psicu^{\mu-1,d}(\Mbar,\Lambda) @>>> \Psicu^{\mu,d}(\Mbar,\Lambda) @>>(\cusym,\cusymb)> {}^{\cu}\Sigma(\Mbar,\Lambda) @>>> 0
\end{CD}
$$
is topologically split exact.

In addition to the cusp-principal symbols, the operators $A(\lambda) = \bigl(A_{i,j}(\lambda)\bigr) \in \Psicu^{\mu,d}(\Mbar,\Lambda)$
have a \emph{conormal symbol} (or normal operator) associated with each singular hypersurface $H \subset \singbdry\Mbar$:
\begin{enumerate}[i)]
\item[iii)] Let $H \subset \singbdry\Mbar$, and assume that $H\cap\regbdry\Mbar \neq \emptyset$. Then the restriction
of $A(\lambda)$ to the collar neighborhood $[0,\eps)\times H \cong \Mbar$ induces an operator
$$
\dot{C}^{\infty}_0\biggl([0,\eps),\begin{array}{c} \dot{C}^{\infty}(H,E) \\ \oplus \\ \dot{C}^{\infty}(\regbdry H,J_-) \end{array}\biggr) \to
\dot{C}^{\infty}\biggl([0,\eps),\begin{array}{c} \dot{C}^{\infty}(H,F) \\ \oplus \\ \dot{C}^{\infty}(\regbdry H,J_+) \end{array}\biggr)
$$
that belongs to the cusp calculus $\Psicu^{\mu,d}([0,\eps)\times H,\Lambda)$ of pseudodifferential boundary value problems
near $H$ considered in Section \ref{subsec-SingularCollar}.
Hence, by Definition \ref{localpseudoconormal}, this operator has a conormal symbol
$N_H(A)(\xi,\lambda) \in \Psicu^{\mu,d}(H,\R\times\Lambda)$, which by definition is the conormal symbol of $A(\lambda)$
with respect to the singular hypersurface $H$.

Recall that the conormal symbol $N_H(A)(\xi,\lambda)$ is a family of boundary value problems \eqref{conormalsymbolcinf}
in the cusp calculus on the manifold with corners $H$.
\item[iv)] Correspondingly, if $H \subset \singbdry\Mbar$ such that $H\cap\regbdry\Mbar = \emptyset$, then the restriction
of the operator $A_{1,1}(\lambda)$ to the collar neighborhood $[0,\eps)\times H$ induces an operator
$$
\dot{C}^{\infty}_0\bigl([0,\eps), \dot{C}^{\infty}(H,E)\bigr) \to
\dot{C}^{\infty}\bigl([0,\eps), \dot{C}^{\infty}(H,F)\bigr)
$$
that belongs to the ordinary cusp calculus $\Psicu^{\mu}([0,\eps)\times H,\Lambda)$ near $H$, see the notes at the
beginning of Section \ref{subsec-SingularCollar}. Analogously to Definition \ref{localpseudoconormal}, this operator
thus has a conormal symbol $N_H(A_{1,1})(\xi,\lambda) \in \Psicu^{\mu}(H,\R\times\Lambda)$.
By definition, we let $N_H(A):= N_H(A_{1,1})$ be the conormal symbol of $A(\lambda)$ with respect to the hypersurface $H$,
a family of cusp pseudodifferential operators
\begin{equation}\label{conormalsymbcinf}
N_H(A)(\xi,\lambda) : \dot{C}^{\infty}(H,E) \to \dot{C}^{\infty}(H,F).
\end{equation}
\end{enumerate}
\end{definition}

\begin{theorem}\label{cuspcomposition}
Let $A_j(\lambda) \in \Psicu^{\mu_j,d_j}(\Mbar,\Lambda)$, $j=1,2$, and assume that the vector bundles fit together such
that the composition $A_1(\lambda)A_2(\lambda)$ is well defined.

Then $A_1(\lambda)A_2(\lambda) \in \Psicu^{\mu_1+\mu_2,\max\{d_1+\mu_2,d_2\}}(\Mbar,\Lambda)$, and the cusp-principal
and conormal symbols of the composition are given by the relations
\begin{equation}\label{princsymbolcomposition}
\left.\begin{aligned}
\cusym(A_1A_2) &= \cusym(A_1)\cusym(A_2), \\
\cusymb(A_1A_2) &= \cusymb(A_1)\cusymb(A_2), \\
N_H(A_1A_2) &= N_H(A_1)N_H(A_2).
\end{aligned}\right\}
\end{equation}
Moreover, for every $\alpha \in \N_0^{N-\ell}$, the class $\x_{\sing}^{\alpha}\Psicu^{\ast,\ast}(\Mbar,\Lambda)$
is a two-sided ideal in the algebra $\Psicu^{\ast,\ast}(\Mbar,\Lambda)$, i.e. whenever any of the $A_j(\lambda)$ above
belongs to this smaller class, so does the composition. Recall that $\singbdry\Mbar$ consists
of $N-\ell$ hypersurfaces.
\end{theorem}
\begin{proof}
Every point $p \in \Mbar$ has an open neighborhood $U(p) \subset \Mbar$ such that, if $p$ is of codimension
$k \in \N$ with $p \in H_{i_1}\cap\ldots\cap H_{i_k}$, where the $H_{i_j} \subset \partial\Mbar$ are $k$ distinct hypersurfaces,
then
$$
U(p) \subset \bigcap\limits_{j=1}^k\bigl([0,\eps)\times H_{i_j}\bigr) \cong \Mbar,
$$
and $U(p)\cap H = \emptyset$ for all hypersurfaces $H \subset \partial\Mbar$ with $H \neq H_{i_j}$, $j=1,\ldots,k$. Moreover,
if $p$ has codimension zero, then $U(p) \cap \partial\Mbar = \emptyset$.
Let $\Mbar = \bigcup\limits_{j=1}^T U(p_j)$ be a finite covering of $\Mbar$ by such neighborhoods, and let
$\{\varphi_j \st j=1,\ldots,T\}$ be a subordinated partition of unity. Choose functions $\varphi_j \prec \psi_j \in
C^{\infty}(\Mbar)$ with $\supp \psi_j \Subset U(p_j)$, $j = 1,\ldots,T$.

To begin with, observe that $\omega A_j(\lambda) \tilde{\omega} \in \Psicu^{-\infty,d_j}(\Mbar,\Lambda)$ for all
functions $\omega,\tilde{\omega} \in C^{\infty}(\Mbar)$ with disjoint supports because
$$
\omega A_j(\lambda) \tilde{\omega} = \sum\limits_{k,l=1}^T \varphi_k\omega A_j(\lambda) \tilde{\omega}\varphi_l,
$$
and every single summand $\varphi_k\omega A_j(\lambda) \tilde{\omega}\varphi_l$ either belongs to
$\Psicu^{-\infty,d_j}_{\infty}(\Mbar,\Lambda)$ by i) of Definition \ref{cuspalgebradef}, or the supports
of the functions $\varphi_k\omega$ and $\tilde{\omega}\varphi_l$ are both contained in a collar neighborhood
$[0,\eps)\times H$ for some $H \subset \singbdry\Mbar$ by construction of the partition of unity.
In the latter case, $\varphi_k\omega A_j(\lambda) \tilde{\omega}\varphi_l$ is in $\Psicu^{-\infty,d_j}(\Mbar,\Lambda)$
by Proposition \ref{localcomposition}.
Recall that the multiplication operators with functions are to be understood according to the convention \eqref{Multop}.

Making use of the defining mapping properties of the residual Green operators of type zero from Section \ref{subsec-ResidualGreen}
and standard arguments in Boutet de Monvel's calculus, we see that if any of the
$A_j(\lambda)$ belongs to $\Psicu_{\infty}^{-\infty,\ast}(\Mbar,\Lambda)$, then so does the composition
$A_1(\lambda)A_2(\lambda)$.

Let us now consider the general composition. We may write $A_1(\lambda)A_2(\lambda)$ as
\begin{equation}\label{comppartition}
A_1(\lambda)A_2(\lambda) = \sum\limits_{j,k,l=1}^T \varphi_j A_1(\lambda) \varphi_k A_2(\lambda) \varphi_l.
\end{equation}
Let us analyze every single summand $\varphi_j A_1(\lambda) \varphi_k A_2(\lambda) \varphi_l$ in
\eqref{comppartition} separately:

Assume that any of the functions $\varphi_j$, $\varphi_k$, or $\varphi_l$ -- say $\varphi_j$ -- is supported
in $\Mbar\setminus\singbdry\Mbar$. Let $\tilde{\psi}_j \in C^{\infty}(\Mbar)$ with $\supp\tilde{\psi}_j \cap \singbdry\Mbar = \emptyset$
and $\psi_j \prec \tilde{\psi}_j$. Write
\begin{align*}
\varphi_j A_1(\lambda) \varphi_k A_2(\lambda) \varphi_l &= \varphi_j A_1(\lambda) \psi_j\varphi_k A_2(\lambda) \varphi_l
+ \bigl(\varphi_j A_1(\lambda) (1-\psi_j)\varphi_k\bigr) A_2(\lambda) \varphi_l \\
&\equiv \varphi_j A_1(\lambda) \psi_j\varphi_k A_2(\lambda) \tilde{\psi}_j\varphi_l
+ \varphi_j A_1(\lambda) \bigl(\psi_j\varphi_k A_2(\lambda) (1-\tilde{\psi}_j)\varphi_l\bigr) \\
&\equiv \varphi_j A_1(\lambda) \psi_j\varphi_k A_2(\lambda) \tilde{\psi}_j\varphi_l.
\end{align*}
Here $\equiv$ means equivalence modulo $\Psicu_{\infty}^{-\infty,\max\{d_1+\mu_2,d_2\}}(\Mbar,\Lambda)$. Observe that
$$
\varphi_j A_1(\lambda) (1-\psi_j), \;
\psi_j A_2(\lambda) (1-\tilde{\psi}_j) \in
\Psicu_{\infty}^{-\infty,\ast}(\Mbar,\Lambda)
$$
by i) of Definition \ref{cuspalgebradef}, and $\Psicu_{\infty}^{-\infty,\ast}(\Mbar,\Lambda)$ is a two-sided ideal.
The functions $\varphi_j$, $\psi_j\varphi_k$, and $\tilde{\psi}_j\varphi_l$ are all supported in
$\Mbar \setminus \singbdry\Mbar$, and consequently the composition $\varphi_j A_1(\lambda) \psi_j\varphi_k A_2(\lambda) \tilde{\psi}_j\varphi_l$
belongs to $\B^{\mu_1+\mu_2,\max\{d_1+\mu_2,d_2\}}(\overline{N},\Lambda)$ by ii) of
Definition \ref{cuspalgebradef} and the composition theorem in Boutet de Monvel's calculus on $\overline{N}$.

The argument for $\varphi_k$ or $\varphi_l$ supported in $\Mbar\setminus\singbdry\Mbar$ is similar, and so
\begin{equation}\label{phiresult}
\varphi_j A_1(\lambda) \varphi_k A_2(\lambda) \varphi_l \in \Psicu^{\mu_1+\mu_2,\max\{d_1+\mu_2,d_2\}}(\Mbar,\Lambda)
\end{equation}
whenever any of the functions $\varphi_j$, $\varphi_k$, or $\varphi_l$ is supported in $\Mbar\setminus\singbdry\Mbar$.

Next assume that the supports of all the functions $\varphi_j$, $\varphi_k$, and $\varphi_l$ have nontrivial intersection
with $\singbdry\Mbar$. If there exists a hypersurface $H \subset \singbdry\Mbar$ with
$$
\supp \varphi_j, \; \supp \varphi_k, \; \supp \varphi_l \subset [0,\eps)\times H \cong \Mbar,
$$
then the composition
$$
\bigl(\varphi_j A_1(\lambda) \varphi_k\bigr)\bigl(\psi_k A_2(\lambda) \varphi_l \bigr) \in \Psicu^{\mu_1+\mu_2,\max\{d_1+\mu_2,d_2\}}([0,\eps)\times H,\Lambda)
$$
by iii), iv) of Definition \ref{cuspalgebradef} and Proposition \ref{localcomposition}, and so \eqref{phiresult} holds
in this case.

By construction of the partition of unity, it remains to consider the case that there exists no $H \subset \singbdry\Mbar$
having nontrivial intersection with the supports of all three functions $\varphi_j$, $\varphi_k$, and $\varphi_l$.
Write
$$
\varphi_j A_1(\lambda) \varphi_k A_2(\lambda) \varphi_l = \varphi_j A_1(\lambda) \varphi_k\psi_l A_2(\lambda) \varphi_l
+ \varphi_j A_1(\lambda) \bigl(\varphi_k(1-\psi_l) A_2(\lambda) \varphi_l\bigr).
$$
We have $(1-\psi_l) A_2(\lambda) \varphi_l \in \Psicu^{-\infty,\ast}(\Mbar,\Lambda)$. Every $\alpha \in \N_0^{N-\ell}$
can be written in the form $\alpha = \alpha_j + \alpha_k + \alpha_l$ with $\alpha_j,\alpha_k,\alpha_l \in \N_0^{N-\ell}$
such that
$$
\tilde{\varphi}_j = \x_{\sing}^{-\alpha_j}\varphi_j, \; \tilde{\varphi}_k = \x_{\sing}^{-\alpha_k}\varphi_k, \text{ and } \tilde{\varphi}_l = \x_{\sing}^{-\alpha_l}\varphi_l \in C^{\infty}(\Mbar).
$$
Now
$$
\varphi_j A_1(\lambda) \varphi_k(1-\psi_l) A_2(\lambda) \varphi_l =
\tilde{\varphi}_j\x_{\sing}^{\alpha_j} A_1(\lambda) \tilde{\varphi}_k\x_{\sing}^{\alpha_k}(1-\psi_l) A_2(\lambda) \tilde{\varphi}_l\x_{\sing}^{\alpha_l},
$$
and consequently this operator is smoothing in the scale of weighted cusp Sobolev spaces, and the range consists of $\dot{C}^{\infty}$-functions
(more precisely, we have to make use of an expansion into operators of type zero and powers of $\partial_+$, and argue for each summand
separately, see Section \ref{subsec-ResidualGreen}).
Thus $\varphi_j A_1(\lambda) \varphi_k(1-\psi_l) A_2(\lambda) \varphi_l \in \Psicu^{-\infty,\ast}_{\infty}(\Mbar,\Lambda)$.

Next pick a function $\tilde{\psi}_l \in C^{\infty}(\Mbar)$ with $\supp \tilde{\psi}_l \Subset U(p_l)$ and $\psi_l \prec
\tilde{\psi}_l$, and write
$$
\varphi_j A_1(\lambda) \varphi_k\psi_l A_2(\lambda) \varphi_l = \varphi_j\tilde{\psi}_l A_1(\lambda) \varphi_k\psi_l A_2(\lambda) \varphi_l
+ \varphi_j(1-\tilde{\psi}_l) A_1(\lambda) \varphi_k\psi_l A_2(\lambda) \varphi_l.
$$
Similar arguments as above give $\varphi_j(1-\tilde{\psi}_l) A_1(\lambda) \varphi_k\psi_l A_2(\lambda) \varphi_l \in
\Psicu^{-\infty,\ast}_{\infty}(\Mbar,\Lambda)$, and, because all functions
$\varphi_j\tilde{\psi}_l$, $\varphi_k\psi_l$, and $\varphi_l$ are supported in one collar neighborhood
$[0,\eps)\times H \cong \Mbar$ for some $H \subset \singbdry\Mbar$, the operator
$$
\varphi_j\tilde{\psi}_l A_1(\lambda) \varphi_k\psi_l A_2(\lambda) \varphi_l \in
\Psicu^{\mu_1+\mu_2,\max\{d_1+\mu_2,d_2\}}([0,\eps)\times H,\Lambda).
$$
Summing up, we have proved \eqref{phiresult} for all possible cases of $\varphi_j$, $\varphi_k$, and $\varphi_l$,
and so the composition theorem is proved. Following the lines of this proof and using Proposition \ref{localcomposition},
we see that the classes $\x_{\sing}^{\alpha}\Psicu^{\ast,\ast}(\Mbar,\Lambda)$ form two-sided ideals for
all $\alpha \in \N_0^{N-\ell}$, and the identities \eqref{princsymbolcomposition} hold.
\end{proof}

\begin{definition}\label{pseudoelliptic}
A boundary value problem $A(\lambda) \in \Psicu^{\mu,d}(\Mbar,\Lambda)$ is called \emph{cusp-elliptic}, if both
the cusp-principal symbol $\cusym(A)$ and the cusp-principal boundary symbol $\cusymb(A)$, see Definition \ref{principalsymbolcusppseudo},
are pointwise bijective.

Moreover, we call $A(\lambda)$ \emph{elliptic}, if in addition all conormal symbols $N_H(A)(\xi,\lambda)$ with
respect to all singular hypersurfaces $H \subset \singbdry\Mbar$ are invertible, i.e. the operator families
\eqref{conormalsymbolcinf} or \eqref{conormalsymbcinf} are bijective for all $(\xi,\lambda) \in \R\times\Lambda$,
respectively.

Recall that $\Lambda \subset \R^q$ is the closure of some open conical subset of $\R^q$, or $\Lambda = \{0\}$.
It would be more precise to reserve the notion of cusp-ellipticity or ellipticity for the case $\Lambda = \{0\}$,
and to call $A(\lambda)$ cusp-elliptic with parameter or elliptic with parameter otherwise.
\end{definition}

\begin{theorem}\label{Parametrix}
Let $A(\lambda) \in \Psicu^{\mu,d}(\Mbar,\Lambda)$ be cusp-elliptic. Then there exists a parametrix
$B(\lambda) \in \Psicu^{-\mu,(d-\mu)_+}(\Mbar,\Lambda)$, $(d-\mu)_+ = \max\{0,d-\mu\}$, such that
$$
\left.\begin{aligned}
A(\lambda)B(\lambda) &- 1 \\
B(\lambda)A(\lambda) &- 1
\end{aligned}\right\}
\in \Psicu^{-\infty,\ast}(\Mbar,\Lambda),
$$
and the types of these remainders are given by the type formula from Theorem~\ref{cuspcomposition}.

Moreover, if $A(\lambda)$ is elliptic, then there exists $R(\lambda) \in \Psicu^{-\infty,(d-\mu)_+}(\Mbar,\Lambda)$
such that
$$
\left.\begin{aligned}
A(\lambda)\bigl(B(\lambda)+R(\lambda)\bigr) &- 1 \\
\bigl(B(\lambda)+R(\lambda)\bigr)A(\lambda) &- 1
\end{aligned}\right\}
\in \Psicu_{\infty}^{-\infty,\ast}(\Mbar,\Lambda),
$$
the space of residual Green operators introduced in Section \ref{subsec-ResidualGreen}.
\end{theorem}
\begin{proof}
Consider the covering
$$
\Mbar = \overline{N} \cup \bigcup\limits_{H\subset \singbdry\Mbar} \bigl([0,\eps)\times H\bigr)
$$
of $\Mbar$, where $\overline{N} = \Mbar\setminus\singbdry\Mbar$ is the regular part of $\Mbar$.
Choose a subordinated partition of unity $\varphi_{\textup{reg}}$, $\varphi_{H}$, $H \subset \singbdry\Mbar$,
and functions $\varphi_{\textup{reg}} \prec \psi_{\textup{reg}}$, $\varphi_{H} \prec \psi_{H}$ that are
compactly supported in $\overline{N}$ or $[0,\eps)\times H$, respectively.

As $A(\lambda)$ is cusp-elliptic, the restriction of $A(\lambda)$ to $\overline{N}$ is elliptic in
Boutet de Monvel's calculus on $\overline{N}$, and, for each singular hypersurface $H \subset \singbdry\Mbar$,
the restriction of $A(\lambda)$ to $[0,\eps)\times H$ is cusp-elliptic in the sense of Definition \ref{localcuspelliptic}.
Choose a parametrix $B_{\textup{reg}}(\lambda)$ in Boutet de Monvel's calculus on $\overline{N}$, as well
as parametrices $B_H(\lambda) \in \Psicu^{-\mu,(d-\mu)_+}([0,\eps)\times H,\Lambda)$ according to Proposition
\ref{localcuspparametrix}, and define
$$
B(\lambda) = \varphi_{\textup{reg}}B_{\textup{reg}}(\lambda) \psi_{\textup{reg}} +
\sum\limits_{H \subset \singbdry\Mbar} \varphi_H B_H(\lambda) \psi_H \in \Psicu^{-\mu,(d-\mu)_+}(\Mbar,\Lambda).
$$
If $H \cap \regbdry\Mbar = \emptyset$, the term $\varphi_H B_H(\lambda) \psi_H$ is by convention a matrix
filled by zeros outside the upper left corner (we employ this convention also further below).
By Theorem \ref{cuspcomposition}, this choice of $B(\lambda)$ furnishes a parametrix of $A(\lambda)$ up to
remainders in the class $\Psicu^{-\infty,\ast}(\Mbar,\Lambda)$ as desired.

Now assume that, in addition, the conormal symbols $N_H(A)(\xi,\lambda)$ are invertible for all $H \subset \singbdry\Mbar$.
Let $H_1,\ldots,H_{N-\ell}$ be an enumeration of the singular boundary hypersurfaces, and let $\x_j$ be the defining
function associated with $H_j$.

We proceed by induction to show that for each $K = 0,\ldots,N-\ell$ there exists a parametrix $B_K(\lambda) \in
\Psicu^{-\mu,(d-\mu)_+}(\Mbar,\Lambda)$ of $A(\lambda)$ such that
\begin{equation}\label{RestK}
\left.\begin{aligned}
A(\lambda)B_K(\lambda) &- 1 \\
B_K(\lambda)A(\lambda) &- 1
\end{aligned}\right\}
\in \bigcap\limits_{\alpha \in T_K} \x_{\sing}^{\alpha}\Psicu^{-\infty,\ast}(\Mbar,\Lambda),
\end{equation}
where
$$
T_K = \{\alpha = (\alpha_1,\ldots,\alpha_{N-\ell}) \in \N_0^{N-\ell} \st \alpha_j = 0 \textup{ for } j > K\}.
$$
The case $K=0$ is fulfilled with $B_0(\lambda):= B(\lambda)$. Assume that we found $B_K(\lambda)$ for some $K < N-\ell$.
We have $A(\lambda)B_K(\lambda) = 1 + R_K(\lambda)$ for some remainder $R_K(\lambda)$ as specified by \eqref{RestK}.
Consequently,
$$
N_{H_{K+1}}(A)^{-1} = N_{H_{K+1}}(B_K) - N_{H_{K+1}}(A)^{-1}N_{H_{K+1}}(R_K),
$$
where
$$
r'(\xi,\lambda):= N_{H_{K+1}}(A)^{-1}N_{H_{K+1}}(R_K) \in \bigcap\limits_{\alpha \in T_K} \x_{\sing}^{\alpha}\Psicu^{-\infty,\ast}(H_{K+1},\R_{\xi}\times\Lambda),
$$
see Theorem \ref{cuspcomposition} and Theorem \ref{largeparaminverse}. Note that these theorems hold by our inductive
approach towards the cusp calculus in view of $\codim H < \codim \Mbar$.

Define an operator $R'(\lambda)$ in the local cusp calculus on $[0,\eps)\times H_{K+1}$ according to the quantization
\eqref{cuspquantization} with the operator-valued symbol $r'(\xi,\lambda)$, and let
$$
B_K'(\lambda):= B_K(\lambda) - \omega(\x_{K+1})R'(\lambda)\tilde{\omega}(\x_{K+1}) \in
\Psicu^{-\mu,(d-\mu)_+}(\Mbar,\Lambda),
$$
where $\omega,\tilde{\omega} \in C_0^{\infty}([0,\eps))$ are cut-off functions that are supported sufficiently close
to the origin. Then $A(\lambda)B_K'(\lambda) = 1 + R'_K(\lambda)$ for some $R'_K(\lambda)$ as specified by
\eqref{RestK}, where in addition $N_{H_{K+1}}(R'_K) \equiv 0$.
Consider the operator
$$
1 + \omega(\x_{K+1})R'_K(\lambda)\tilde{\omega}(\x_{K+1}) \in 1 + \bigcap\limits_{\alpha \in T_K}
\x_{\sing}^{\alpha}\Psicu^{-\infty,\ast}([0,\eps)\times H_{K+1},\Lambda).
$$
According to Proposition \ref{localcuspparametrix}, there exists
$$
\tilde{R}_K(\lambda) \in \bigcap\limits_{\alpha \in T_K}
\x_{\sing}^{\alpha}\Psicu^{-\infty,\ast}([0,\eps)\times H_{K+1},\Lambda)
$$
such that $\bigl(1 + \omega(\x_{K+1})R'_K(\lambda)\tilde{\omega}(\x_{K+1})\bigr)\bigl(1 + \tilde{R}_K(\lambda)\bigr)$
belongs to
$$
1 + \bigcap\limits_{\substack{\alpha \in T_K \\ j \in \N_0}}
\x_{K+1}^{j}\x_{\sing}^{\alpha}\Psicu^{-\infty,\ast}([0,\eps)\times H_{K+1},\Lambda).
$$
Define
$$
B_{K+1}(\lambda):= B'_K(\lambda)\bigl(1 + \omega(\x_{K+1})\tilde{R}_K(\lambda)\tilde{\omega}(\x_{K+1})\bigr) \in
\Psicu^{-\mu,(d-\mu)_+}(\Mbar,\Lambda).
$$
Then
$$
A(\lambda)B_{K+1}(\lambda) \in 1 + \bigcap\limits_{\alpha \in T_{K+1}} \x_{\sing}^{\alpha}\Psicu^{-\infty,\ast}(\Mbar,\Lambda)
$$
by construction. The same arguments apply for the left parametrix, and the induction is therefore complete.
Hence
$$
\left.\begin{aligned}
A(\lambda)B_{N-\ell}(\lambda) &- 1 \\
B_{N-\ell}(\lambda)A(\lambda) &- 1
\end{aligned}\right\}
\in \bigcap\limits_{\alpha \in \N_0^{N-\ell}} \x_{\sing}^{\alpha}\Psicu^{-\infty,\ast}(\Mbar,\Lambda) =
\Psicu^{-\infty,\ast}_{\infty}(\Mbar,\Lambda),
$$
and the theorem is proved.
\end{proof}

\begin{corollary}\label{elliptFredh}
Let $A \in \Psicu^{\mu,d}(\Mbar)$ be elliptic. Then the extension
\begin{equation}\label{Areal}
A : \begin{array}{c} \x_{\sing}^{\alpha}\Hcu^{s}(\Mbar,E) \\ \oplus \\ \x_{\sing}^{\alpha}\Hcu^{s}(\regbdry\Mbar,J_-) \end{array} \to
\begin{array}{c} \x_{\sing}^{\alpha}\Hcu^{s-\mu}(\Mbar,F) \\ \oplus \\ \x_{\sing}^{\alpha}\Hcu^{s-\mu}(\regbdry\Mbar,J_+) \end{array}
\end{equation}
is a Fredholm operator for all $s > \max\{\mu,d\} - \frac{1}{2}$ and $\alpha \in \R^{N-\ell}$.
\end{corollary}
\begin{proof}
Let $B + R \in \Psicu^{-\mu,(d-\mu)_+}(\Mbar)$ be a parametrix of $A$ up to residual Green operators. By
\eqref{Sobspaceembedding} the residual Green operators are compact in the cusp Sobolev spaces. Hence $B+R$
is an inverse of \eqref{Areal} up to compact operators.
\end{proof}

\begin{theorem}\label{largeparaminverse}
Assume that $\Lambda \neq \{0\}$, and let $A(\lambda) \in \Psicu^{\mu,d}(\Mbar,\Lambda)$ be cusp-elliptic with
parameter $\lambda \in \Lambda$. Then, for all $\lambda \in \Lambda$ outside possibly a compact set
$K \subset \Lambda$, the operator $A(\lambda)$ is invertible in the cusp Sobolev spaces \eqref{extensionSobCont}
for $s > \max\{\mu,d\}-\frac{1}{2}$ and all $\alpha \in \R^{N-\ell}$, as well as in the spaces \eqref{initialcusp}
of smooth functions that vanish to infinite order at $\singbdry\Mbar$.

Moreover, for any open neighborhood $K \subset U(K)$, there exists a parameter-dependent parametrix
$B(\lambda) \in \Psicu^{-\mu,(d-\mu)_+}(\Mbar,\Lambda)$ of $A(\lambda)$ such that $A(\lambda)B(\lambda) = \id$ and
$B(\lambda)A(\lambda) = \id$ for $\lambda \notin U(K)$, i.e. $B(\lambda) = A(\lambda)^{-1}$ both in
the spaces \eqref{extensionSobCont} and \eqref{initialcusp}.
In particular, if $A(\lambda)$ happens to be invertible for all $\lambda \in \Lambda$, then $A(\lambda)^{-1} \in
\Psicu^{-\mu,(d-\mu)_+}(\Mbar,\Lambda)$.
\end{theorem}
\begin{proof}
Let $B'(\lambda) \in \Psicu^{-\mu,(d-\mu)_+}(\Mbar,\Lambda)$ be a parameter-dependent parametrix of $A(\lambda)$
according to Theorem \ref{Parametrix}, i.e.
$$
\left.\begin{aligned}
A(\lambda)B'(\lambda) - 1 &= R'_r(\lambda) \\
B'(\lambda)A(\lambda) - 1 &= R'_l(\lambda)
\end{aligned}\right\}
\in \Psicu^{-\infty,\ast}(\Mbar,\Lambda) = \S\bigl(\Lambda,\Psicu^{-\infty,\ast}(\Mbar)\bigr).
$$
Consequently, $A(\lambda)B'(\lambda)$ and $B'(\lambda)A(\lambda)$ are invertible for large $\lambda \in \Lambda$, i.e.
$A(\lambda)$ is invertible for all $\lambda \in \Lambda$ outside possibly some compact set $K \subset \Lambda$.

Let $\chi \in C^{\infty}(\Lambda)$ with $\chi \equiv 0$ in some neighborhood of $K$, and $\chi \equiv 1$ outside
the given neighborhood $U(K)$ of $K$. Define
$$
B(\lambda) = B'(\lambda) - B'(\lambda)R'_r(\lambda) + R'_l(\lambda)\chi(\lambda)A(\lambda)^{-1}R'_r(\lambda).
$$
As both $R'_l(\lambda)$ and $R'_r(\lambda)$ are regularizing, we deduce that
$$
R'_l(\lambda)\chi(\lambda)A(\lambda)^{-1}R'_r(\lambda) \in \Psicu^{-\infty,(d-\mu)_+}(\Mbar,\Lambda),
$$
and so $B(\lambda) \in \Psicu^{-\mu,(d-\mu)_+}(\Mbar,\Lambda)$. Moreover,
$A(\lambda)B(\lambda) = \id$ for $\lambda \notin U(K)$ by construction. This proves the theorem.
\end{proof}



\begin{thebibliography}{10}

\bibitem{AgmonDouglisNirenberg}
S.~Agmon, A.~Douglis, and L.~Nirenberg, \emph{Estimates near the boundary for solutions of elliptic
partial differential equations satisfying general boundary conditions I, II}, Commun.~Pure Appl.~Math.
\textbf{12} (1959), 623--727, and \textbf{17} (1964), 35--92.

\bibitem{AV63}
M.~Agranovich and M. Vishik, \emph{Elliptic problems with a
parameter and parabolic problems of general type}, Russ. Math. Surveys
\textbf{19} (1963), 53--159.

\bibitem{AmLauNis}
B.~Ammann, R.~Lauter, and V.~Nistor, \emph{Pseudodifferential operators on manifolds
with a Lie structure at infinity}, preprint 2003 (math.AP/0304044 at arXiv.org).

\bibitem{BdM71}
L.~Boutet de Monvel, \emph{Boundary problems for pseudo-differential operators},
Acta Math. \textbf{126} (1971), 11--51.

\bibitem{Cordes}
H.~Cordes, \emph{A global parametrix for pseudo-differential operators over ${\mathbb R}^n$, with applications}, SFB 72, preprint 90, Bonn, 1976.

\bibitem{Dauge}
M.~Dauge, \emph{Elliptic boundary value problems on corner domains}, vol.~1341 of \emph{Lecture
Notes in Mathematics}. Springer-Verlag, Berlin, 1988.

\bibitem{Eskin}
G.~Eskin, \emph{Boundary value problems for elliptic pseudodifferential equations}, Translations of Mathematical Monographs, vol.~52. American Mathematical Society, Providence, R.I., 1981.

\bibitem{Fedosov}
B.~Fedosov, B.-W.~Schulze, and N.~Tarkhanov, \emph{The index of elliptic operators on manifolds with conical
points}, Sel.~Math. New Ser. \textbf{5} (1999), 467--506.

\bibitem{GKM1}
J.~Gil, T.~Krainer, and G.~Mendoza, \emph{Geometry and spectra of 
closed extensions of elliptic cone operators}, preprint 2004 (math.AP/0410178 at arXiv.org), to appear in
Canadian J.~Math.

\bibitem{GKM2}
J.~Gil, T.~Krainer, and G.~Mendoza, \emph{Resolvents of elliptic
cone operators}, preprint 2004 (math.AP/0410176 at arXiv.org).

\bibitem{GrubbBuch} 
G.~Grubb, \emph{Functional calculus of pseudodifferential boundary problems},
2nd ed., Progress in Mathematics, vol.~65. Birkh{\"a}user, Basel, 1996.

\bibitem{GruSe95}
G.~Grubb and R.~Seeley, \emph{Weakly parametric pseudodifferential operators
and Atiyah-Patodi-Singer boundary problems}, Invent. Math. \textbf{121} (1995),
no. 3, 481--529.

\bibitem{JerKen}
D.~Jerison and C.~Kenig, \emph{The inhomogeneous Dirichlet problem in Lipschitz
domains}, J. Funct. Anal. \textbf{130} (1995), no.~1, 161--219.

\bibitem{KaSchu03}
D.~Kapanadze and B.-W.~Schulze, \emph{Crack Theory and Edge Singularities},
Mathematics and its Applications, vol.~561, Kluwer Academic Publishers,
Dordrecht-Boston-London, 2003.

\bibitem{K1} 
V.~Kondratyev, \emph{Boundary problems for
elliptic equations in domains with conical or angular points},
Trans. Mosc. Math. Soc. \textbf{16} (1967), 227--313.

\bibitem{KozMazRos}
V.~Kozlov, V.~Mazya, and J.~Ro{\ss}mann, \emph{Spectral problems associated with corner
singularities of solutions to elliptic equations}, vol.~85 of \emph{Mathematical Surveys and
Monographs}. American Mathematical Society, Providence, RI, 2001.

\bibitem{KraiParabBVP}
T.~Krainer, \emph{On the inverse of parabolic boundary value problems for large times}, Japanese J. Math. \textbf{30} (2004), 91--163.

\bibitem{KraConicResolvents}
T.~Krainer, \emph{Resolvents of elliptic boundary problems on conic manifolds}, preprint 2005 (math.AP/0503021 at arXiv.org).

\bibitem{LauterMoroianu}
R.~Lauter and S.~Moroianu, \emph{The index of cusp operators on manifolds with corners},
Ann. Global Anal. Geom. \textbf{21}, no. 1, (2002), 31--49.

\bibitem{Lesch}
M.~Lesch, \emph{Operators of Fuchs type, conical singularities, and asymptotic methods}, Teubner Texte zur Mathematik,
vol.~136, Teubner-Verlag, Leipzig, 1997.

\bibitem{LioMag}
J.-L.~Lions and E.~Magenes, \emph{Non-homogeneous boundary value problems and applications}, Grundlehren der mathematischen Wissenschaften, vol.~181--183. Springer-Verlag, New York-Heidelberg, 1972, 1973.

\bibitem{Loya}
P.~Loya, \emph{The index of $b$-pseudodifferential operators on manifolds with corners}, Ann.~Global Anal.~Geom.
\textbf{27} (2005), 101--133.

\bibitem{LoyaPark04}
P.~Loya and J.~Park, \emph{Boundary problems for Dirac-type operators on manifolds with
multi-cylindrical end boundaries}, preprint 2004.

\bibitem{MazNazPlam}
V.~Mazya, S.~Nazarov, and B.~Plamenevsky, \emph{Asymptotic theory of elliptic boundary value
problems in singularly perturbed domains {I}, {II}}, vol.~111 and 112 of \emph{Operator Theory: Advances
and Applications}, Birkh{\"a}user Verlag, Basel, 2000.

\bibitem{MazzeoMelrose}
R.~Mazzeo and R.~Melrose, \emph{Pseudodifferential operators on manifolds with fibred boundaries}, Asian J.~Math.
\textbf{2} (1998), 833--866.

\bibitem{RBM2} 
R.~Melrose, \emph{The Atiyah-Patodi-Singer index
theorem}, Research Notes in Mathematics, A~K~Peters, Ltd.,
Wellesley, MA, 1993.

\bibitem{MM} 
R.~Melrose and G.~Mendoza, \emph{Elliptic operators
of totally characteristic type}, MSRI Preprint, 1983.

\bibitem{MelroseNistor}
R.~Melrose and V.~Nistor, \emph{Homology of pseudodifferential operators I. Manifolds with boundary},
to appear in Amer.~J.~Math.

\bibitem{MimiTay}
D.~Mitrea, M.~Mitrea, and M.~Taylor, \emph{Layer potentials, the Hodge Laplacian, and global boundary
problems in nonsmooth Riemannian manifolds}, vol.~713 of \emph{Memoirs of the American Mathematical
Society}. American Mathematical Society, Providence, RI, 2001.

\bibitem{MiNistor}
M.~Mitrea and V.~Nistor, \emph{Boundary value problems and layer potentials on manifolds with
cylindrical ends}, preprint 2004 (math.AP/0410186 at arXiv.org).

\bibitem{MiTayVas}
M.~Mitrea, M.~Taylor, and A.~Vasy, \emph{Lipschitz domains, domains with corners, and the
Hodge Laplacian}, preprint 2004 (math.AP/0408438 at arXiv.org).

\bibitem{Mueller}
W.~M{\"u}ller, \emph{On the $L^2$-index of Dirac operators on manifolds with corners of codimension two, I},
J.~Differential Geom. \textbf{44} (1996), 97--177.

\bibitem{Nistor}
V.~Nistor, \emph{Singular integral operators on non-compact manifolds and analysis on polyhedral
domains}, preprint 2004 (math.AP/0402322 at arXiv.org).

\bibitem{Parenti}
C.~Parenti, \emph{Operatori pseudodifferenziali in ${\mathbb R}^n$ e applicazioni}, Annali Mat. Pura ed Appl. \textbf{93} (1972), 359--389.

\bibitem{Schrohe99}
E.~Schrohe, \emph{Fr{\'e}chet algebra techniques for boundary value problems on noncompact
manifolds: Fredholm criteria and functional calculus via spectral invariance},
Math.~Nachr. \textbf{199} (1999), 145--185.

\bibitem{SchroSchu}
E.~Schrohe and B.-W.~Schulze, \emph{Boundary value problems in {B}outet de {M}onvel's algebra for manifolds with conical singularities {I},{II}},
Math. Top., vol.~5 and 8, Akademie Verlag, Berlin, (1994, 1995), pp.~97--209, 70--205.

\bibitem{SzNorthHolland} 
B.-W.~Schulze, \emph{Pseudo-differential operators on manifolds with
singularities}, Studies in Mathematics and its Applications, vol.~24. 
North-Holland Publishing Co., Amsterdam, 1991.

\bibitem{SzWiley98}
\bysame, \emph{Boundary value problems and singular
pseudo-differential operators}, Pure and Applied Mathematics.
John Wiley \& Sons, Ltd., Chichester, 1998.

\bibitem{SchulzeKyoto}
\bysame, \emph{Operators with symbol hierarchies and iterated asymptotics}, Publications
of the RIMS, Kyoto University, \textbf{38} (2002), 735--802.

\bibitem{Seeley}
R.~Seeley, \emph{The resolvent of an elliptic boundary problem},
Amer. J. Math. \textbf{91} (1969), 889--920. 

\bibitem{Shubin}
M.~Shubin, \emph{Pseudodifferential operators and spectral theory}, Springer Verlag, Princeton, N.J., 1987.

\bibitem{Vaillant}
B.~Vaillant, \emph{Index and spectral theory for manifolds with generalized fibred cusps},
Ph.D. thesis, Universit{\"a}t Bonn, 2001.

\end{thebibliography}
\end{document}